\documentclass[10pt,journal,compsoc]{IEEEtran}
\usepackage{lineno,hyperref}
\usepackage{graphicx,times,amsmath}
\usepackage{setspace}
\usepackage{algorithm}
\usepackage{algorithmic}
\usepackage{enumerate}
\usepackage{subfigure}
\usepackage{stmaryrd}
\usepackage{amsfonts}
\usepackage{multirow}
\usepackage{fancyvrb}
\usepackage{array,tabularx,,supertabular}
\usepackage{amssymb}
\usepackage{amsmath,amssymb}
\usepackage{cases}
\usepackage{ulem,booktabs}
\usepackage{color}
\usepackage{subfigure}
\usepackage{caption}
\usepackage{bm}
\usepackage{float}
\usepackage{multicol}
\usepackage{epstopdf}
\usepackage{algorithm, algorithmic}
\usepackage{setspace}
\usepackage{makecell}
\usepackage{ulem}
\ifCLASSOPTIONcompsoc
  \usepackage[nocompress]{cite}
\else
  \usepackage{cite}
\fi
\ifCLASSINFOpdf
\else
\fi
\begin{document}
%
\title{A Multiobjective State Transition Algorithm Based on Decomposition}
%
%
%
%

\author{Xiaojun~Zhou,\
        Yuan~Gao,
        Shengxiang Yang,\
        Chunhua Yang,
        and~ Jiajia Zhou
\IEEEcompsocitemizethanks{\IEEEcompsocthanksitem
X. J. Zhou, Y. Gao, C. H. Yang and J. J. Zhou are with the School of Automation, Central South University, Changsha 410083, China, and X. J. Zhou is also with the State Key Laboratory of Synthetical Automation for Process Industries, Shenyang 110000, China. \protect\\
\IEEEcompsocthanksitem
S. X. Yang is with the School of Computer Science and Informatics, De Montfort University, Leicester LE1 9BH, U.K. \protect\\
 \protect\\
}
\thanks{Manuscript received XX, 2021; revised XX, 2021.}}

%
%

\markboth{Journal of \LaTeX\ Class Files,~Vol.~14, No.~8, August~2015}%
{Shell \MakeLowercase{\textit{et al.}}: Bare Demo of IEEEtran.cls for Computer Society Journals}
%



\IEEEtitleabstractindextext{%
\begin{abstract}
Aggregation functions largely determine the convergence and diversity performance of multi-objective evolutionary algorithms in decomposition methods. Nevertheless, the traditional Tchebycheff function does not consider the matching relationship between the weight vectors and candidate solutions. In this paper, the concept of matching degree is proposed which employs vectorial angles between weight vectors and candidate solutions. Based on the matching degree, a new modified  Tchebycheff aggregation function is proposed, which integrates matching degree into the Tchebycheff aggregation function. Moreover, the proposed decomposition method  has the same functionality with the Tchebycheff aggregation function. Based on the proposed decomposition approach, a new multiobjective optimization algorithm named  decomposition based multi-objective state transition algorithm  is proposed. Relevant experimental results show that the proposed algorithm is highly competitive in comparison with other state-of-the-art multiobjetive optimization algorithms.
\end{abstract}

\begin{IEEEkeywords}
Multi-objective optimization, decomposition, evolutionary algorithms, matching degree, Tchebycheff approach, state transition algorithm
\end{IEEEkeywords}}

\maketitle

\IEEEdisplaynontitleabstractindextext

%
\IEEEpeerreviewmaketitle

\IEEEraisesectionheading{\section{Introduction}\label{sec:introduction}}

%
%
%
%
\IEEEPARstart {M}any real-world engineering problems often involve the optimization of several different conflicting objectives \cite{nag2014asmiga}.
They are often referred to as  multi-objective optimization problems (MOPs). It is necessary to find optimization approaches to solve these problems  effectively and obtain solutions with trade-offs  among different objectives. Evolutionary algorithms are well suited for solving MOPs by obtaining a solution set in one single run. Over the past two decades, multi-objective evolutionary algorithms (MOEAs) have developed rapidly \cite{yen2003dynamic}.

According to the selection strategy, MOEAs can be divided into three categories \cite{zhou2011multiobjective}: (i) MOEAs based on Pareto dominance, (ii) MOEAs based on decomposition, (iii) MOEAs based on indicators. Dominance-based MOEAs have been prevalent in recent decades
\cite{zitzler1999multiobjective} such as NSGA \cite{srinivas1994muiltiobjective}, PAES \cite{knowles1999pareto}, SPEA-II \cite{zitzler2001spea2}, and NSGA-II \cite{deb2002fast}. However, the domination principle will be too weak to provide an adequate selection pressure. Many methods based on performance indicators are proposed, e.g.,  hypervolume (HV) indicator \cite{li2012achieving}, the S-metric selection-based evolutionary multiobjective algorithm \cite{beume2007sms}, and HypE \cite{bader2011hype}. A  indicator-based evolutionary algorithm (IBEA) was proposed that can be combined with arbitrary indicators \cite{zitzler2004indicator}. In the IBEA, there is no need for any  diversity preservation mechanism such as fitness sharing.
Unfortunately, this type of MOEAs tends to be time-consuming when calculating performance indicators
\cite{while2006faster}.

MOEAs based on decomposition obtain increasing attention in recent years \cite{LIANG201850,ZHU2018263,QI2017844}.  In MOEA/D \cite{zhang2007moea}, a MOP is decomposed into a set of subproblem, and each solution is associated with a subproblem.
With the growing complexity of real-world MOPs whose Pareto Fronts (PFs) tend to be irregular, some weight vector adjustment strategies are proposed to enhance the existing MOEA/D on solving such MOPs \cite{ma2020survey}.  Several variants of MOEA/D are proposed to enhance the selection strategy for each subproblem \cite{tan2013moea,sato2014inverted,li2014interrelationship,asafuddoula2014decomposition,gee2014online}.
MOEA/D-M2M works by dividing PF into a set of segments and solving them separately \cite{liu2013decomposition}. RVEA is guided by a set of predefined reference vectors \cite{cheng2016reference}.
MOEA/DD takes advantage of dominance- and decomposition-based approaches, it can be able to  balance convergence and diversity \cite{li2014evolutionary}. A decomposition-based idea is employed in  NSGA-III to maintain population diversity, and the concept of Pareto dominance is adopted to maintain the population convergence \cite{deb2013evolutionary}.

Aggregation functions largely determine the convergence and diversity performance of MOEAs in decomposition methods. Many decomposition approaches are proposed to make algorithms get better convergence and distribution. The Inverted PBI (IPBI) decomposition method can better approximate widely spread Pareto fronts \cite{sato2015analysis}. Adaptive penalty scheme (APS) and subproblem-based penalty scheme (SPS) are  proposed  to solve the problem that PBI needs to adjust parameters appropriately, and they can improve algorithm convergence and diversity \cite{yang2017improving}. Moreover, a Tchebycheff decomposition-based MOEA with $l_2$-norm is proposed. In MOEA/D-MR, both the ideal points and the nadir points are adopted in decomposition methods to obtain Pareto optimal solutions \cite{wang2017use}. Meanwhile, there are also plenty of achievements regarding to other improvements on decomposition based algorithms. However, the traditional Tchebycheff method does not consider the matching relationship between the candidate solutions and weight vectors, which may cause the Pareto optimal solution obtained not uniformly distributed and make better solutions not be retained.

In addition to decomposition approach, search ability is vital for MOEAs. Lots of heuristic algorithms with different search  strategies have been proposed in recent years. For example, MOEA/D-DE \cite{li2009multiobjective}, multi-objective particle swarm optimizatio (MOPSO) algorithm  \cite{zapotecas2011multi}, etc. State transition algorithm (STA) \cite{zhou2012state} is proposed and presents excellent performance compared with other global optimization algorithms \cite{han2017new,han2018two,huang2018hybrid,zhou2018statistical,zhou2019dynamic,zhou2020hybrid,zhou2020kernel,zhou2020using,zhou2020feature}.  When the objective functions and their PFs are non-convex, the state transformation operators proposed in STA are advantageous for exploration and exploitation. Various state transformation operators can be used for global search, local search and heuristic search. Alternative use of local search and global search, which can quickly converge to the PFs for saving the search time.

In this paper, the influence between
the matching of  weight vectors and candidate solutions on the
update process of candidate solutions is analyzed in detail.
The concept of matching degree is proposed, and vectorial angle is employed to evaluate the matching degree between weight vectors and candidate solutions. Based on the matching degree, a modified Tchebycheff approach is proposed. Furthermore, a decomposition based multi-objective state transition algorithm (MOSTA/D) is proposed based on the modified Tchebycheff approach. The
main new contributions of this paper can be summarized as
follows.

\noindent
\hangafter=1
\setlength{\hangindent}{2em}
1)  State transformation operators are adopted to reproduce candidate solutions in a collaborative manner. These state transformation operators not only can be controlled with search region but also can balance the global search and local search.

\noindent
\hangafter=1
\setlength{\hangindent}{2em}
2) The concept of matching degree is proposed which considers the matching relationship between weight vectors and candidate solutions. Based on the matching degree, a new decomposition approach named modified Tchebycheff approach is proposed.

\noindent
\hangafter=1
\setlength{\hangindent}{2em}
3) A new decomposition based multi-objective state transition algorithm is proposed. Verified by several benchmark test functions, the proposed algorithm is valid and effective to solve MOPs with complex Pareto set shapes, and it can obtain Pareto optimal solutions with good convergence and diversity.

The structure of this paper is organized as follows. Section 2
introduces some background knowledge and analyzes the matching relationship between
candidate solutions and weight vectors. The details of the
proposed MOSTA/D are described in Section 3. Section 4
presents experiments and discussions. Finally, the conclusion and future work are given in Section 5.

%

\section{Related work}
In this section, a brief review of basic definitions and Tchebycheff decomposition are presented firstly. Then, an introduction to the analysis of the decomposition aggregation function in this paper is given, including the matching relationship  how to influence the updating process.
\subsection{Basic definitions}
The MOP solved in this paper can be formulated as follows:
\begin{eqnarray}
\min\limits_{\bm x \in {\Omega}}& \bm F(\bm x)=( f_1(\bm x),f_2(\bm x),\cdots,f_m(\bm x))^T
\label{eq:2.2_multi}
\end{eqnarray}
where $\Omega$ is the decision space and $\bm{x} = (x_1,\cdots,x_n) \in \Omega$ is an $n$-dimensional decision variable vector which represents a solution to the target MOP. $\bm{F}(\bm{x}) :\Omega \rightarrow \mathbb{R}^{m}$ denotes the $m$-dimensional objective vector of the solution $\bm{x}$. \\

\subsection{Tchebycheff decomposition}
The scalar optimization problem of Tchebycheff approach is described in the following:
\begin{equation}
\label{che}
\min {\rm{ }}{g^{{\rm{te}}}}({\bm{x}}|{\bm{\lambda }},{{\bm{z}}^{\bm{*}}}) = \mathop {\max }\limits_{1 \le i \le m} \left\{ {{\lambda _i}|{f_i}({\bm{x}}) - z_i^*|} \right\}\\
\end{equation}
where $\bm{\lambda}$ = $(\lambda_1, \lambda_2,\cdots,\lambda_m)^T$ is a weight vector whose length is equal to the number of the objective function, with $\lambda_i \geq 0 $ and $\Sigma_{i=1}^{m}${$\lambda_i$} =1 for all $i$ = $1,\cdots,m$. $z^{\ast}=  (z_{ 1}^{\ast},\dots , z_{m}^{\ast})^T$ is the reference point. There exists a correspondence that each Pareto optimal solution is an optimal solution of objective function Eq. (\ref{che}).

\begin{figure}[!htbp ]
  \centering
  \includegraphics[scale=0.45,trim=150 100 0 100]{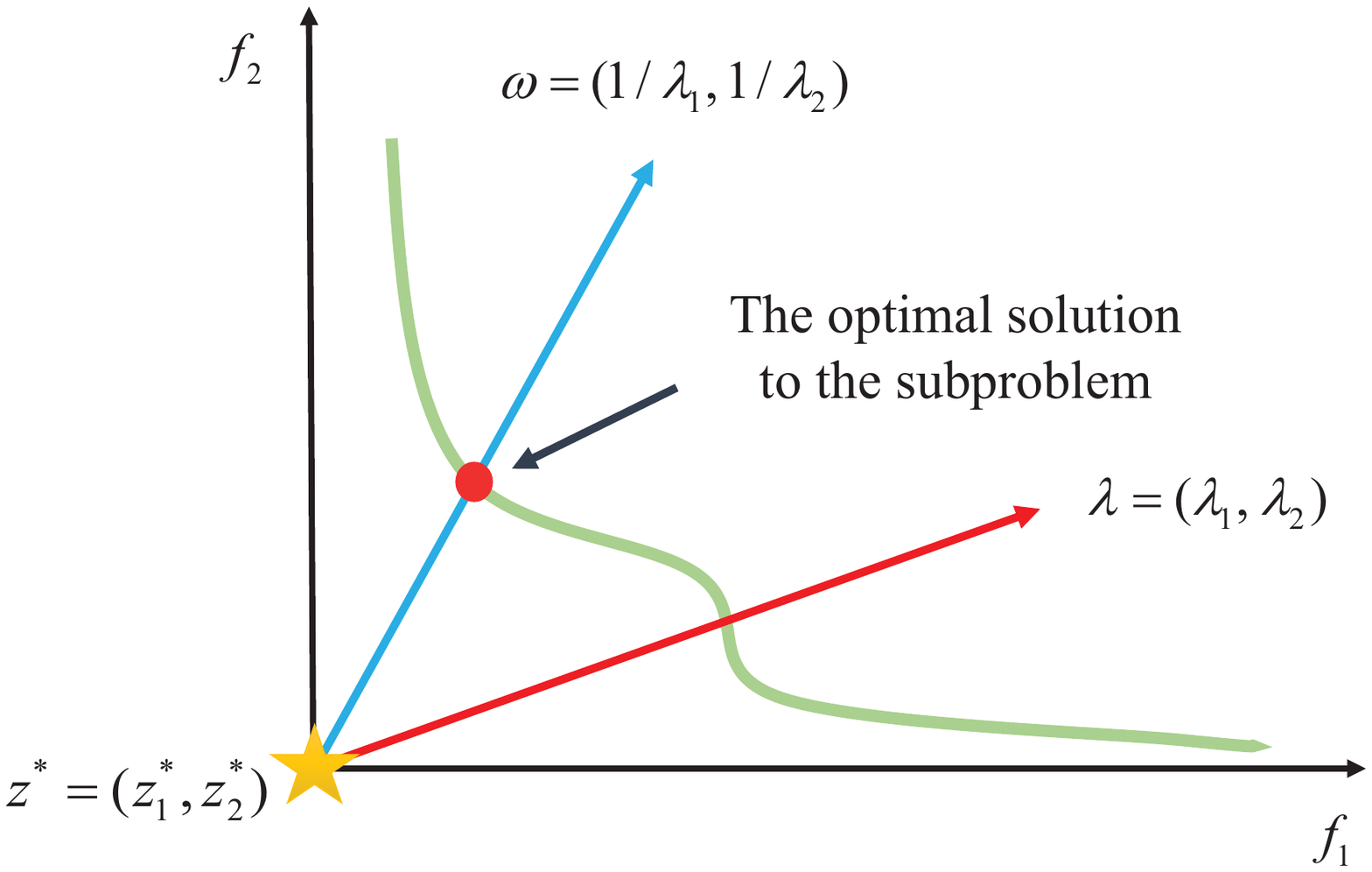}
 \caption{The schematic diagram of the Tchebycheff decomposition approach}
 \label{sdgp}
 \end{figure}

The above theorem is expressed intuitively in Fig. \ref{sdgp} based
on a two-objective MOP. The green curve is the PF
of the problem. The optimal solution in the objective space of
the scalar subproblem with weight vector $\lambda$ is collinear with $\omega$.

\subsection{Motivations}
When the ideal reference point is fixed, the decomposition aggregation function can be regarded as a function of weight vectors and candidate solutions.
Hence, the decomposition aggregation function value of one candidate solution
with different weight vectors are quite different.  Weighted sum aggregation function value is the projection on the weight vector, which owns good geometric performance and is easy to understand. Therefore, the weighted sum approach is tasken as an example to give further analysis of influence of decomposition aggregation function in updating process.

\begin{figure}[!htbp]
  \centering
  \includegraphics[scale=0.7,trim=250 200 100 170]{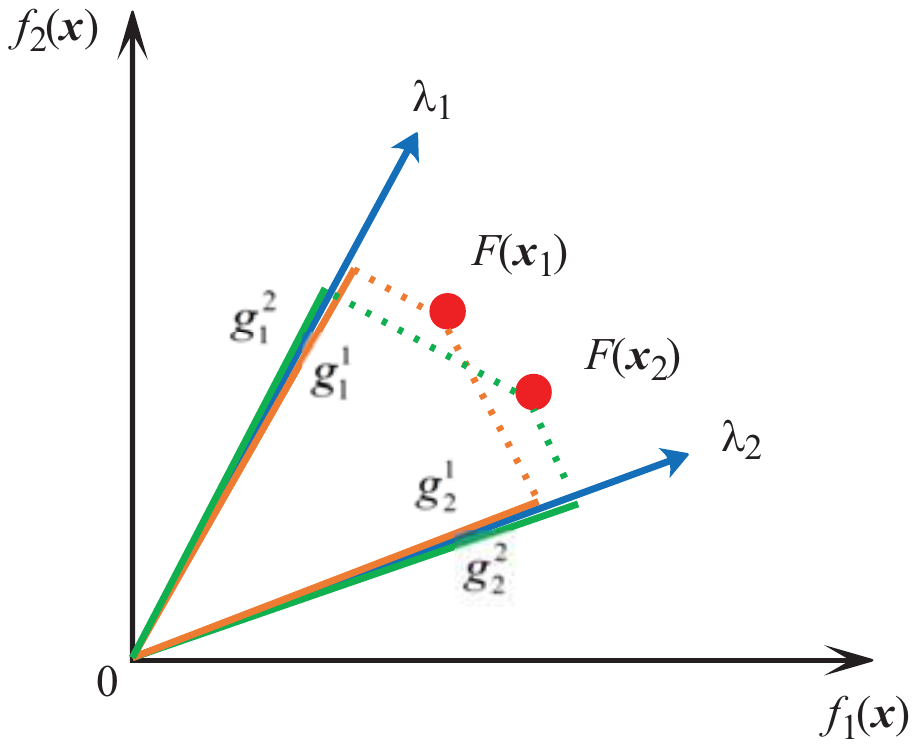}
\caption{The schematic diagram of the matching relationship}
 \label{Motivation}
\end{figure}

The matching relationship between weight vectors and candidate solutions demonstrated
in Fig.~\ref{Motivation}.
$\bm{x_1}$  is one of the solutions in
the primitive population and $\bm{x_2}$ is one of the solutions in the
new population.
$\bm{f(\bm{x_1})}$ and $\bm{f(x_2)}$ are the objective functions.
$\bm{\lambda_1}$ and $\bm{\lambda_2}$ are different weight vectors.
The objective function of the scalar optimization
problem based on the weighted sum approach of two solutions
is considered. $g_i^j$ represents the objective
function value of the scalar optimization problem of  $\bm{x_j}$ when
it matches the weight vector $\bm{\lambda_i}$.

From Fig.~\ref{Motivation}, we can find that if  $\bm{x_1}$ and $\bm{x_2}$ match with  weight vector $\bm{\lambda_1}$, then $g_1^1$ is bigger than $g_1^2$, so according to the replacement criteria of aggregation function, $\bm{x_1}$ would be replaced by $\bm{x_2}$. However, if  $\bm{x_1}$ and $\bm{x_2}$ match with the weight vector $\bm{\lambda_2}$, then $g_2^1$ is bigger than $g_2^2$, so according to the
replacement criteria of aggregation function, $\bm{x_1}$ would not be replaced by $\bm{x_2}$. It can be found that the results of updating population process are different when the same candidate solution is matched with different weight vectors.

Therefore, the matching degree between
weight vectors and candidate solutions  should be carefully considered in the updating procedure with the decomposition approach. Moreover, the
decomposition approach should comprehensively consider the
optimal matching degree. Considering that the Tchebycheff aggregation function is the most commonly used decomposition
approach, in this paper, the matching degree in the Tchebycheff aggregation function
will be analyzed in detail.

\section{Decomposition based multi-objective state transition algorithm}
The details of the proposed MOSTA/D are described in this section.  A
new decomposition approach based on the matching degree,
named modified Tchebycheff decomposition. The updating process not only updates
the population, but also strengthens the information communication among potential excellent solutions.  Further explanation will be illustrated in the following.

\subsection{Initialization of weight vectors}
\begin{figure}[!htbp]
  \centering
  \includegraphics[width=7cm]{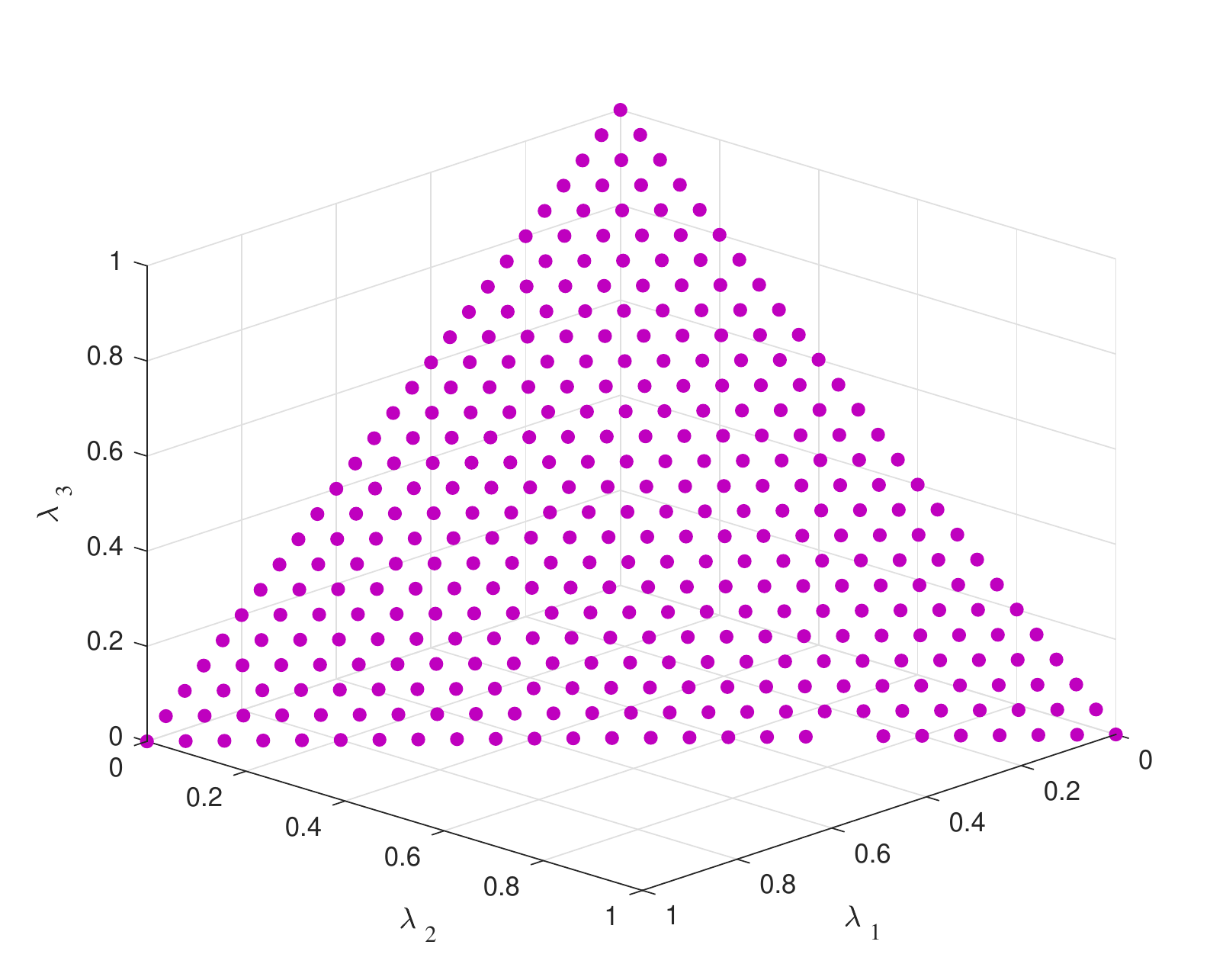}
\caption{The schematic diagram of weight vectors}
   \label{L}
 \end{figure}

For each weight vector $\bm{\lambda}^i$ = ($\lambda^i_1, \lambda^i_2,\cdots,\lambda^i_m$), the elements takes values from $\{0/N, 1/N, 2/N,\cdots,N/N\}$ under the condition of $\sum_{i=1}^m{\lambda^i_i} = 1$. For a MOP with $m$ objectives, $N= C_{N+m-1}^{m-1}$ is the number of such vectors,
where $N$ is a user-defined positive integer. If $N$ is larger than
the number, we can sample weight vectors up to the
number. The neighborhood $B$ between subproblems
are obtained by calculating Euclidean distances.
The weight vector obtained by uniform sampling is shown in Fig. \ref{L}.

For each subproblem, it has $T$ subproblems in its neighborhood. The initial population $P= \{{\bm{x}^{(1)}},{\bm{x}^{(2)}}, \cdots ,{\bm{x}^{(n)}}\}$ is randomly sampled. A candidate
solution $\bm{x}^i$ is assigned to a subproblem randomly. The ideal reference point $\bm{z}^*$ is set as    $\bm{z}^* = \min \{f_i(\bm{x})|\bm{x} \in P\}$.

\subsection{Reproduction}
Four state transformation operators are used for generating candidate solutions. Different state transformation operators can be used for global search, local search and heuristic search. Alternative use of different operators, so that the state transition algorithm can quickly find the global optimal solution with a certain probability.
\begin{itemize}
\item[$\bullet$] Rotation transformation
\begin{equation}
s_{k+1}=  s_{k}+ \alpha \frac{1}{\omega \| s_{k}\|_{2}} R_{r}  s_{k},
\end{equation}
where $s_{k} \in \rm{R}^n$ is a candidate solution, $\alpha$ is a positive constant, called the rotation factor,
$R_{r}$ is a random variable
with its components obeying the uniform distribution in the range of [0,1],
and $\|\cdot\|_{2}$ is the 2-norm of a vector.
\item[$\bullet$]  Translation transformation
\begin{equation}
s_{k+1} = s_{k}+  \beta  R_{t}  \frac{ s_{k}-  s_{k-1}}{\| s_{k}-  s_{k-1}\|_{2}},
\end{equation}
where $\beta$ is a positive constant, called the translation factor.
$R_{t}$ is a random variable with its components obeying the uniform distribution in the range of [0,1].
\item[$\bullet$]
Expansion transformation
\begin{equation}
s_{k+1}=  s_{k}+\gamma R_{e}s_{k},\\
\end{equation}
where $\gamma$ is a positive constant, called the expansion factor. $R_{e}$ is a random diagonal matrix with its entries
obeying the Gaussian distribution.
\item[$\bullet$]  Axesion transformation\\
\begin{equation}
s_{k+1} =  s_{k}+  \delta  R_{a}   s_{k},\\
\end{equation}
where $ \delta$ is a positive constant, called the axesion factor
$R_{a}$ is a random diagonal matrix whose entries obey a Gaussian distribution with variable variance and only one random position has a nonzero value.
\end{itemize}

\subsection{A new decomposition approach based on  Tchebycheff and matching degree}
The candidate solutions matching different weight vectors
would cause the difference of decomposition aggregation
function values. Furthermore, it affects the selection and updating
process of candidate solutions. Hence, the matching degree of
weight vectors and candidate solutions is critical for selection
in decomposition based algorithms. However, the Tchebycheff decomposition approach not explicitly highlights the matching degree between
the weight vectors and candidate solutions, which may cause the Pareto optimal solution not uniformly distributed. In this section, a new decomposition approach
based on matching degree is proposed, which comprehensively
takes into account the Tchebycheff decomposition approach
and the relationship. It can be demonstrated that the proposed approach
has the same functionality with the Tchebycheff decomposition approach.

First, the following lemma explains the geometric properties
of the Tchebycheff decomposition approach:\\
\textbf{Lemma 1} \textit{It is assumed that the target Pareto front of the multiobjective problems to be solved is piecewise continuous.
If the straight line $\bm{L_1}$ :
$\frac{f_1-z_1^*}{\frac{1}{\lambda_1}}$ = $\frac{f_2-z_2^*}{\frac{1}{\lambda_2}}$= $\cdots$ = $\frac{f_m-z_m^*}{\frac{1}{\lambda_m}}$
, taking $f_1, f_2, \cdots, f_m$ as variables,
has an intersection with the PF, then the intersection point is the optimal solution \bm{$x^*$ } to the scalar
subproblem with weight vector $\lambda$ = $(\lambda_1, \lambda_2, \cdots, \lambda_m)$. $z^*$ =
$(z^*_1, z^*_2, \cdots, z^*_m)$ is the ideal reference vector of optimization problems. }


From the theorem mentioned above, it can be concluded that $\bm{F(x^*|\bm{\lambda})} - \bm{z}^*$ is collinear with $\bm{\omega}$ = $(1/\lambda_1, 1/\lambda_2,\cdots,1/\lambda_m)$,
where $\bm{F(x^*|\bm{\lambda})} $ is the objective function values vector. Here, cosine value of vectorial angle is introduced to represent the relationship between $\bm{F(x^*|\bm{\lambda})} $ and $\bm{\omega}$.\\
\textbf{Definition 5} (\textit{Cosine Value of Vectorial Angle})
Cosine value of vectorial angle of $\bm{v}$ and $\bm{u}$ is defined as follows:
\begin{equation}
\cos(\bm{v},\bm{u})  \triangleq (\frac{\bm{v}\cdot{\bm{u}}}{||\bm{v}|| \cdot ||\bm{u}||})
\end{equation}

If $\bm{v}$ and $\bm{u}$ are collinear, the absolute value of $\cos(\bm{v},\bm{u})$ is equal to 1.
The more consistent the direction between  $\bm{v}$ and $\bm{u}$,  the closer the absolute value of $\cos(\bm{v},\bm{u})$ to 1.
Therefore, focusing on $\bm{F}(x^*|\bm{\lambda})$ and $\bm{\omega}$, the absolute value of $\cos(\bm{F(x^*)}|\bm{\lambda},\bm{\omega})$ is equal to 1. Whereas, if the absolute value of $\cos(\bm{F(x)}|\bm{\lambda},\bm{\omega})$ is not equal to 1, $\bm{x}$ is not the optimal solution.

 \begin{figure}[!htbp]
  \centering
  \includegraphics[scale=0.4,trim=300 100 200 100]{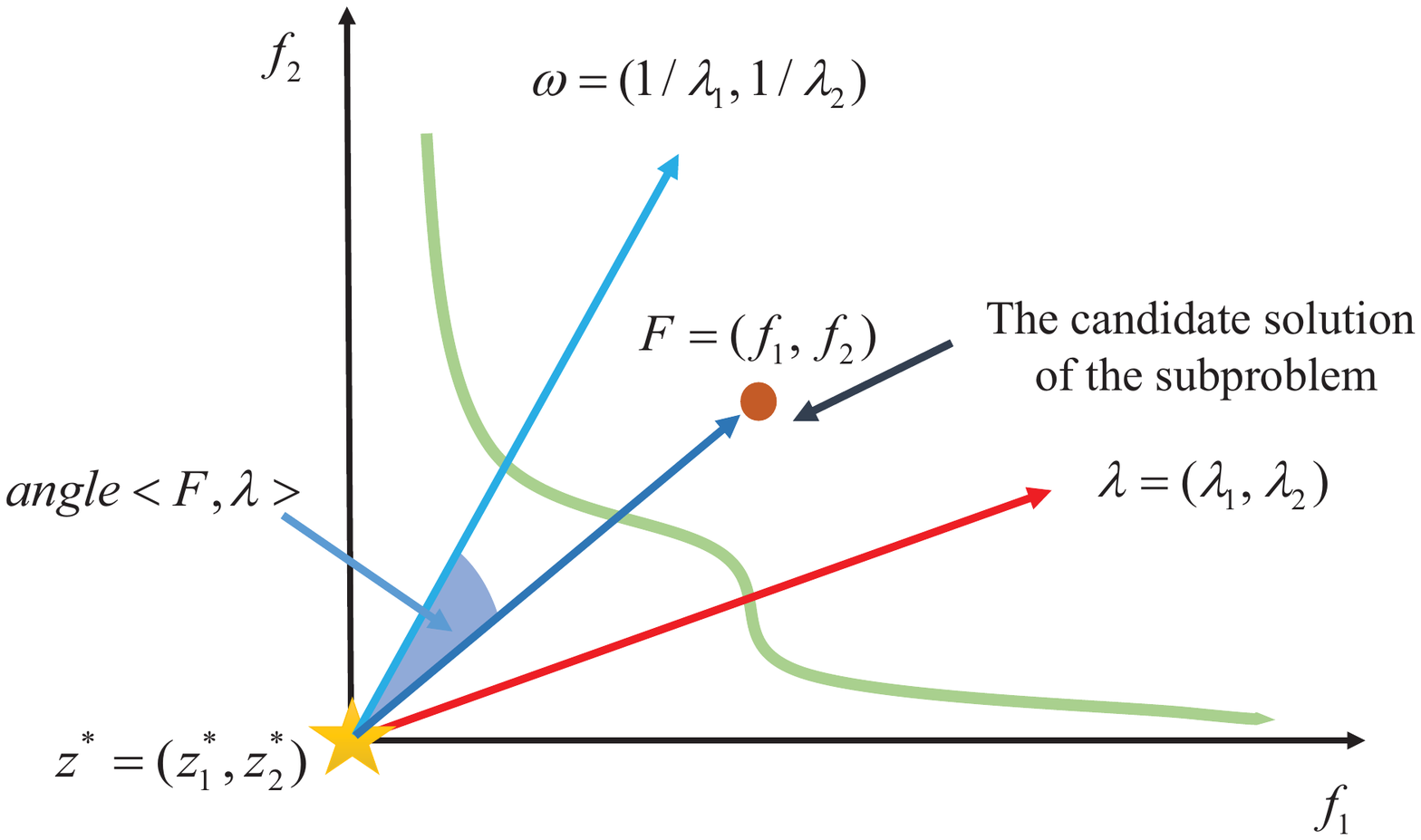}
\caption{The schematic diagram of vectorial angle in decomposition approach}
   \label{angle}
 \end{figure}

Based on the above analysis, the matching degree of weight
vectors and candidate solutions based on the vectorial angle
is defined as follows:\\
\textbf{Definition 6} (\textit{Matching Degree})
The matching degree $\phi$ between the weight vector $\bm{\lambda}$  and candidate solution $\bm{x}$ is:

\begin{equation}
\begin{aligned}
\phi  = & \left |  \cos(\omega,\bm{F(\bm{x}|\bm{\lambda})-\bm{z^*}}) -1 \right |\\
=& \left |  \frac{\bm{\omega}\cdot {(\bm{F(x|\bm{\lambda})-\bm{z^*})}}}{||\bm{\omega}|| \cdot ||\bm{F(x|\bm{\lambda})-\bm{z^*}}||} -1 \right |
\end{aligned}
\end{equation}
where $\bm{\omega}$ = $(1/\lambda_1, 1/\lambda_2,\cdots,1/\lambda_m)^T$, which is a row vector and
$F(\bm{x}|\bm{\lambda})$ = $(f_1(\bm{x}|\bm{\lambda}), f_2(\bm{x}|\bm{\lambda}),\cdots,f_M(\bm{x}|\bm{\lambda}))^T$, which is
the vector of objective function values of the MOP with
solution $\bm{x}$, and $\bm{F}(\bm{x}|\bm{\lambda})$ is also a row vector.
The schematic diagram of the matching degree is shown in Fig. \ref{angle}. $angle< F, \lambda>$ represents the vectorial angle between $F$ and $\lambda$.  The smaller the  $angle<F, \lambda>$, the closer $\phi$ is to 0. In objective function space, $\bm{F}(\bm{x})$  is  closer with $\lambda$ , which
is considered that  candidate solution $x$ is more suitable matching with the weight vector $\lambda$.
On the contrary, The larger the  $angle<F, \lambda>$, $\bm{F}(\bm{x})$  is  farer with $\lambda$, which
is considered that  candidate solution $x$, is less suitable matching with the weight vector $\lambda$.

Based on  Tchebycheff and matching degree, a new decomposition approach is proposed as follows:

\begin{algorithm}
	\renewcommand{\algorithmicrequire}{\textbf{Input:}}
	\renewcommand{\algorithmicensure}{\textbf{Output:}}
	\caption{The pseudo-code of the proposed method}
	\label{alg}
	\begin{algorithmic}[1]
		\REQUIRE
        $MOP$; the popsize of population $N$;  the number of the weight vectors in the neighborhood $Nh$; search enforcement $SE$;  a stopping criterion;
		\ENSURE
        the Pareto optimal solution $P^*$;
    \STATE \textbf{Step:1 Initialization} \\
    Initialization weight vector $\bm{\lambda}$, initial population $P$ and calculate the objective function $FP$,  the ideal reference point $\bm{z}^*$, and neighborhood set $B(i)$ for weight vector $\lambda^i$ corresponding to solution $P(i)$. Initialization  parameters of the state transformation operators $\alpha$, $\gamma$, $\beta$, $\delta$, upper and lower bounds of state transformation factors $\alpha _{max}$, $\alpha _{min}$, $\gamma _{max}$, $\gamma _{min}$, and decay rate of state transformation factors $fc_{\alpha}$, $fc_{\gamma}$.
    \STATE\textbf{Step:2 Reproduction} \\
    Randomly select $T$ candidate solutions from the population $P$ as initial solutions ${P}'$, and each candidate solution generates $SE$ new solutions by one or more state transformation operators and form new population $Q$, and calculate the objective function $F(Q)$.
     \STATE $Q_{1}$ $\leftarrow$ rotation~(${P}'$, $SE$, $\alpha$)
     \STATE $Q_{2}$ $\leftarrow$  expansion~(${P}'$, $SE$, $\gamma$)
     \STATE $Q_{3}$ $\leftarrow$  axesion~(${P}'$, $SE$, $\delta$)
     \STATE $Q \leftarrow  \left \{Q_{1},Q_{2},Q_{3}\right \}$
      \STATE \textbf{Step:3 Update}
        \STATE Generate a random number $\delta$
        \IF  { $\delta < $0.5}
             \STATE Randomly select $Nh$ weight vectors as $B(i)$
        \ENDIF
                 \IF { $\alpha < \alpha _{min}$}
             \STATE $\alpha =\varepsilon_{\alpha max}$
        \ENDIF
             \IF { $ \gamma < \gamma_{min}$}
             \STATE $ \gamma = \gamma _{max}$
        \ENDIF
        \STATE ${\bm{z}^*} \leftarrow \min ({\bm{z}^*},F(Q))$
                \FOR {$i$ = 1 to $N$}
                    \FOR {$j$ $\in$ $B(i)$}
                        \IF {$g^{tmd}$($\bm{P(i)}|\bm{\lambda^j}$, $\bm{z}^*$) $\geq$ $g^{tmd}$($\bm{Q(j)}|\bm{\lambda^j}$, $\bm{z}^*$)}
                               \STATE$State$ $\leftarrow$  translation($\bm{P(i)}$, $\bm{Q(i)}$, $SE$, $\beta$)
                                 \IF   {$g^{tmd}$($\bm{Q(i)}|\bm{\lambda^j}$, $\bm{z}^*$) $\geq$ $g^{tmd}$($State|\bm{\lambda^j}$, $\bm{z}^*$)}
                                   \STATE $\bm{P(i)}$ $\leftarrow$  $State$
                                    \IF   {$g^{tmd}$($\bm{Q(i)}|\bm{\lambda^j}$, $\bm{z}^*$) $\geq$ $g^{tmd}$($State|\bm{\lambda^j}$, $\bm{z}^*$)}
                                     \STATE$\bm{P(i)}$ $\leftarrow$  $\bm{Q(i)}$
                                     \ENDIF
                               \ENDIF
                        \ENDIF
                   \ENDFOR
                \ENDFOR
           \STATE {$\alpha = \alpha /fc_{\alpha}$}
             \STATE {$\varepsilon_\gamma = \varepsilon_\gamma /fc_{\gamma}$}
	        \STATE \textbf{Step:4 Stopping Criteria}
                 \STATE  If stopping criteria is satisfied, then stop and output. Otherwise, go to Steps 2 and 3.
  \end{algorithmic}
\label{alg}
\end{algorithm}

\begin{equation}
\textrm{min} \qquad g^{tmd}(\bm{x}|\bm{\lambda}) = g^{te}(\bm{x}|\bm{\lambda})  \cdot(1+\phi)\\
\label{matching1}
\end{equation}
where $g^{te}(\bm{x}|\bm{\lambda})$ is the Tchebycheff aggregation function defined in general, $\phi$ is the matching degree.

\textit{Proposition 1}: Let  $\bm{x_{te}^*}$ be the optimal solution of Eq. (\ref{che}) and $\bm{x_{tmd}^*}$ be the optimal solution of Eq. (\ref{matching1}) for a fixed $\bm{\lambda}$. Correspondingly, let $g^{te}(\bm{x_{te}^*})$ be  the optimal value of Eq. (\ref{che}) and $g^{tmd}(\bm{x_{tmd}^*}|\bm{\lambda})$ is optimal value of Eq. (\ref{matching1}). It can be concluded that $\bm{x_{tmd}^*}$ is the same as $\bm{x_{te}^*}$ and $g^{tmd}(\bm{x_{tmd}^*}|\bm{\lambda})$ is equal to $g^{te}(\bm{x_{te}^*}|\bm{\lambda})$.

\textit{Proof}: From the construction of $\phi$, when $\bm{\lambda}$ is fixed, it can be found that $\phi\geq 0$. Therefore, it can be concluded
that $g^{tmd}(\bm{x}|\bm{\lambda})\geq g^{th}(\bm{x}|\bm{\lambda})$.
Let $\bm{x_{te}^*}$ be the optimal solution of Eq.(\ref{che}) and $\bm{x_{tmd}^*}$ the optimal solution of Eq.(\ref{matching1}). Therefore, $g^{tmd}(\bm{x}|\bm{\lambda})\geq g^{te}(\bm{x}|\bm{\lambda}) \geq g^{te}(\bm{x_{te}^*}|\bm{\lambda})$. So, $g^{tmd}(\bm{x^*}|\bm{\lambda})$ has a lower bound and has a minimum value. When the $g^{te}(\bm{x}|\bm{\lambda})$ obtains the optimal value, the  value of $\phi$ is equal to 0. Meanwhile,  $g^{tmd}(\bm{x}|\bm{\lambda})= g^{te}(\bm{x}|\bm{\lambda})$ and $g^{tmd}(\bm{x}|\bm{\lambda})$ obtains the optimal value. Hence, $\bm{x_{tmd}^*}= \bm{x_{te}^*}$ and $g^{te}(\bm{x_{te}^*}|\bm{\lambda}) = g^{tmd}(\bm{x_{te}^*}|\bm{\lambda})$.

\subsection{Update procedure and complexity analysis}
As mentioned above, the updating procedure is shown in Algorithm \ref{alg}. Compared with the objective value of the
scalar subproblems, the better solutions are stored in $P$. It is
worth noting that the update process strengthens the information communication of potential excellent solutions. When the
offspring candidate solutions are compared with parent candidate
solutions based on the proposed decomposition approach and
offspring candidate solutions are superior to parent candidate
solutions, those solution are considered as potential excellent
solutions and strengthened search will be activated. Those
solutions will act as initial solutions and will be transformed by
translation transformation operators to generate new solutions.
More information are put in lines 9-34 of Algorithm \ref{alg}.

According to Algorithm  \ref{alg}, the main computational complexity is determined by updating candidate solutions. The population size is $N$,  and the number of the weight vectors in the neighborhood is $Nh$. Besides, when evaluating
the modified Tchebycheff approach, the time complexity of translation transformation is $O(SE)$. In summary, the Update step need $O(N \cdot Nh \cdot SE)$ comparisons. The Initialization step and Reproduction step can finish in linear time.

\section{Experimental analysis}
In this section, several experiments are conducted to verify the convergence and diversity of the proposed algorithm.  The population
size is 200 and each algorithm runs 30 times independently. The stopping criterion
of all algorithms is that the maximum number of
objective function evaluations reaches 10$e^5$. The four state transformation operators adopted in MOSTA/D
are insensitive to the values of parameters in the range of 0.5-0.9. To make a fair
comparison of various algorithms, a compromise parameter value is usually adopted within this range. The parameters of comparing algorithms are set to their default values in PlatEMO \cite{tian2017platemo}. All the benchmark test functions \cite{huband2006review} are shown in Table \ref{Bechmark}.

 The Wilcoxon rank sum test is adopted to  compare the results at a significance level of 0.05.  Symbol ``-" indicates that the compared algorithm  is significantly outperformed by  MOSTA/D, while ``+" means that MOSTA/D  is significantly outperformed by the compared algorithm. Finally, ``$\approx$" means that there is no statistically significant difference between them.

\subsection{Performance metrics}
Two  performance metrics are adopted in assessing the performance of the compared algorithms on benchmark test functions.  The details are given in the following:

1) Modified Inverted Generational Distance ($IGD^+$) Metric \cite{10.1007/978-3-319-15892-1_8}
\begin{equation}
IGD^+(P^*,P) = \frac{\sum_{P\in P^*}{\text{dist}(P,PF)}}{|P^*|}
\end{equation}
where ${\text{dist}(P,PF)}$ denotes the nearest distance from $P$ to the solution $Y$ in $PF$, and the distance  is calculated by $\sqrt{\sum_{j=1}^{m} \max \left(Y_{j}-P_{j}, 0\right)^{2}}$. $|P^*|$ is the number of solutions in $P^*$. Obviously, the smaller the value of $IGD^+$ is, the better convergence and diversity algorithm has.

2) Hypervolume ($HV$) Metric \cite{while2006faster}
\begin{equation}
HV(P,{\mathbf{z}^ * }) = {\rm{volume}}(\bigcup\limits_{{\bf{x}} \in P} {[{f_1}(\mathbf{x}),z_1^ * ] \times } ...[{f_m}(\mathbf{x}),z_m^ * ])
\end{equation}
where ${\rm{volume \left ( \cdot  \right )}}$ indicates the Lebesgue measure. The larger is
the HV value, the better is the quality of $P$ for approximating
the $PF$. $\mathbf{z}^{\ast}=(z_{1}^{\ast}, \ldots, z_{m}^{\ast})^{T}$  is a reference point. In the experiments, $\mathbf{z}^ *$ is set to be 1.2 times the maximum value of the objective function in the PF.

\begin{table*}
\begin{center}
  \caption{Benchmark test functions of multi-objectives optimization algorithms}
  \label{Bechmark}
  \scalebox{0.76}{\begin{tabular}{cccccccc}
    \toprule
  Problem   & dimension (D) & variable domain & objective functions \\
  \midrule
   \multirow{4}*{\makecell{P1}}  & \multirow{4}*{\makecell{7}}& \multirow{4}*{\makecell{[0,1]}}    &  $f_1(x) = 0.5x_{1}x_{2}(1+g(x))  $\\
                                    &                              &                                 &  $f_2(x) = 0.5x_{1}(1-x_{2})(1+g(x))  $\\
                                    &                              &                                 &  $f_3(x) = 0.5(1-x_{1})(1+g(x))  $\\
                                     &                              &                                & $g(x)= 100(5 + \sum\nolimits_{i = 3}^n {{{({x_i} - 0.5)}^2} - \cos (20\pi ({x_i} - 0.5))} )$\\
     \midrule
   \multirow{4}*{\makecell{P2}}  & \multirow{4}*{\makecell{12}}& \multirow{4}*{\makecell{[0,1]}}    &  $f_1(x) = \cos (0.5\pi {x_1})\cos (0.5\pi {x_2})(1 + g(x))$\\
                                    &                              &                                 &  $f_2(x) = \cos (0.5\pi {x_1})\sin (0.5\pi {x_2})(1 + g(x))$\\
                                    &                              &                                 &  $f_3(x) = \sin (0.5\pi {x_2})(1 + g(x))$\\
                                     &                              &                                & $g(x)= \sum\nolimits_{i = 3}^n {{{({x_i} - 0.5)}^2}} $\\

  \midrule
   \multirow{4}*{\makecell{P3}}  & \multirow{4}*{\makecell{12}}& \multirow{4}*{\makecell{[0,1]}}    &  $f_1(x) = \cos (0.5\pi x_1^a)\cos (0.5\pi x_2^a)(1 + g(x))$\\
                                    &                              &                                 &  $f_2(x) = \cos (0.5\pi x_1^a)\sin (0.5\pi x_2^a)(1 + g(x))$\\
                                    &                              &                                 &  $f_3(x) = \sin (0.5\pi x_1^a)(1 + g(x))$\\
                                     &                              &                                & $g(x)= {\sum\nolimits_{i = 3}^n {({x_i} - 0.5)} ^2}$\\

    \midrule
  \multirow{3}*{\makecell{P4}}  & \multirow{3}*{\makecell{13}}& \multirow{3}*{\makecell{$z_{i=1: n, \max }=2 i$}}    &  $t_{i=1: n}^{1} =\rm {s\_multi}\left(y_{i}, 30,10,0.35\right)$\\
                                    &                             &                                 & $t_{i=1: M-1}^{2} =\rm {r\_sum }\left(\left\{y_{(i-1) k /(M-1)+1}, \ldots, y_{i k /(M-1)}\right\},\{1, \ldots, 1\}\right)$\\
                                     &                              &                                & $t_{M}^{2} =\rm{r\_sum }\left(\left\{y_{k+1}, \ldots, y_{n}\right\},\{1, \ldots, 1\}\right)$\\

    \midrule
     \multirow{4}*{\makecell{P5}}  & \multirow{4}*{\makecell{10}}& \multirow{4}*{\makecell{[0,1]}} & $f_1(x) =   (1+g(x))x_1$ \\
                                    &                              &                                &  $f_2(x) =   (1+g(x))(1-x_1^2)$\\
                                    &                              &                                & $g(x)=10*\sin(\pi x_1)*\sum_{i=2}^n{\frac{|t_i|}{1+e^{5|t_i|}}}$\\
                                    &                              &                                & $t_i = x_i - \sin(0.5\pi x_1)$ \\
    \midrule

      \multirow{4}*{\makecell{P6}}  & \multirow{4}*{\makecell{10}}& \multirow{4}*{\makecell{[0,1]}} &$ f_1(x) =(1+g(x))\cos(\frac{\pi x_i}{2})$ \\
                                    &                              &                                &  $f_2(x) = (1+g(x))(1-x_1^2)$ \\
                                    &                              &                                & $g(x)=10*\sin(\pi x_1)*\sum_{i=2}^n{\frac{|t_i|}{1+e^{5|t_i|}}}$ \\
                                    &                              &                                & $t_i = x_i - \sin(0.5\pi x_1)$ \\
    \midrule

          \multirow{4}*{\makecell{P7}}  & \multirow{4}*{\makecell{10}}& \multirow{4}*{\makecell{[0,1]}} &$f_1(x) =(1+g(x))x_1$\\
                                    &                              &                                &  $f_2(x) = (1-x_1^{0.5}{\cos(2\pi x_1)}^2)(1+g(x))$  \\
                                    &                              &                                & $g(x)=10*\sin(\pi x_1)*\sum_{i=2}^n{\frac{|t_i|}{1+e^{5|t_i|}}}$ \\

                                    &                              &                                & $t_i = x_i - \sin(0.5\pi x_1)$ \\
    \midrule

    \multirow{4}*{\makecell{P8}}  & \multirow{4}*{\makecell{30}}& \multirow{4}*{\makecell{[0,1]}} &$f_1(x) =  (1+g(x))(1-x_1)  $\\
                                    &                              &                                &  $f_2(x)=0.5* (1+g(x))(x_1+\sqrt{x_1}\cos^2(4\pi x_1)) $\\
                                    &                              &                                & $g(x)= 2\sin(0.5\pi x_i)(n-1+\sum_{i=2}^n{y_i^2-\cos(2\pi y_i)}) $\\
                                     &                              &                               & $y_{i=2:n}=x_i-\sin(0.5\pi x_i) $\\
    \midrule

    \multirow{4}*{\makecell{P9}}  & \multirow{4}*{\makecell{30}}& \multirow{4}*{\makecell{[0,1]}} &$f_1(x) =  (1+g(x))x_1  $\\
                                    &                              &                                 & $f_2(x)=0.5* (1+g(x))(1-x_1^{0.1}+(1-\sqrt{x_2})^2\cos^2(3\pi x_1)) $\\
                                    &                              &                                 & $g(x)= 2\sin(0.5\pi x_i)(n-1+\sum_{i=2}^n{y_i^2-\cos(2\pi y_i)}) $\\
                                     &                              &                                & $y_{i=2:n}=x_i-\sin(0.5\pi x_i) $\\
    \midrule

    \multirow{4}*{\makecell{P10}}  & \multirow{4}*{\makecell{30}}& \multirow{4}*{\makecell{[1,4]}}    &  $f_1(x) = x_1(1+g(x))/\sqrt{x_2x_3}  $\\
                                    &                              &                                 &  $f_2(x) = x_2(1+g(x))/\sqrt{x_1x_3}  $\\
                                    &                              &                                 &  $f_3(x) = x_3(1+g(x))/\sqrt{x_1x_2}  $\\
                                     &                              &                                & $g(x)=\sum_{i=4}^n{(x_i-2)^2} $\\

    \midrule
   \multirow{4}*{\makecell{P11}}  & \multirow{4}*{\makecell{12}}& \multirow{4}*{\makecell{[0,1]}}    &  $f_1(x) = \cos (0.5\pi {x_1})\cos (0.5\pi {x_2})(1 + g(x))$\\
                                    &                              &                                 &  $f_2(x) = \cos (0.5\pi {x_1})\sin (0.5\pi {x_2})(1 + g(x))$\\
                                    &                              &                                 &  $f_3(x) = \sin (0.5\pi {x_1})(1 + g(x))$\\
                                     &                              &                                & $g(x)= 100(5 + \sum\nolimits_{i = 3}^n {{{({x_i} - 0.5)}^2} - \cos (20\pi ({x_i} - 0.5))} )$\\
   \midrule
  \multirow{6}*{\makecell{P12}}  & \multirow{6}*{\makecell{14}}& \multirow{6}*{\makecell{$z_{i=1: n, \max }=2 i$}}    &  $t_{i=1: k}^{1} =y_{i}$\\
  &                              &                                 &  $t_{i=k+1: n}^{1} = \rm{s\_linear} \left(y_{i}, 0.35\right)$\\
                                    &                              &                                 &  $t_{i=1: k}^{2} =y_{i}$\\
                                     &                              &                                & $t_{i=k+1: k+l / 2}^{2} =\rm {r\_nonsep} \left(\left\{y_{k+2(i-k)-1}, y_{k+2(i-k)}\right\}, 2\right)$\\
                                     &                              &                                & $t_{i=1: M-1}^{3} =\rm{r\_sum}\left(\left\{y_{(i-1) k /(M-1)+1}, \ldots, y_{i k /(M-1)}\right\},\{1, \ldots, 1\}\right)$ \\
                                     &                              &                                &
                                     $t_{M}^{3} =\rm{r\_sum}\left(\left\{y_{k+1}, \ldots, y_{k+l / 2}\right\},\{1, \ldots, 1\}\right)$\\

 \midrule
  \multirow{3}*{\makecell{P13}}  & \multirow{3}*{\makecell{13}}& \multirow{3}*{\makecell{$z_{i=1: n, \max }=2 i$}}    &  $t_{i=1: n}^{1} =\quad \rm{ s\_decept }\left(y_{i}, 0.35,0.001,0.05\right)$\\
&                              &                                 &  $t_{i=1: M-1}^{2} =\rm {r\_sum }\left(\left\{y_{(i-1) k /(M-1)+1}, \ldots, y_{i k /(M-1)}\right\},\{1, \ldots, 1\}\right)$\\
                                     &                              &                                & $t_{M}^{2} =\rm{r\_sum }\left(\left\{y_{k+1}, \ldots, y_{n}\right\},\{1, \ldots, 1\}\right)$\\

  \midrule
   \multirow{4}*{\makecell{P14}}  & \multirow{4}*{\makecell{13}}& \multirow{4}*{\makecell{$z_{i=1: n, \max }=2 i$}}    &  $t_{i=1: k}^{1} =y_{i}$\\
  &                              &                                 &  $t_{i=k+1: n}^{1} = \rm{s\_linear} \left(y_{i}, 0.35\right)$\\
   &                              &                                 & $t_{i=1: M-1}^{2} =\rm { r\_nonsep }\left(\left\{y_{(i-1) k /(M-1)+1}, \ldots, y_{i k /(M-1)}\right\}, k /(M-1)\right)$\\
   &                              &                                 &$t_{M}^{2} =\rm{ r\_nonsep }\left(\left\{y_{k+1}, \ldots, y_{n}\right\}, l\right)$\\
\bottomrule
  \end{tabular}}
  \end{center}
\end{table*}

%
%

\subsection{Validation of the proposed decomposition method}
P1-P4 are adopted to demonstrate the effectiveness of the modified Tchebycheff aggregation function based on matching degree in multi-objective optimization algorithms. Traditional Tchebycheff aggregation function is adopted in  state trasition algorithm to conduct comparision experiments which is called MOSTA/D-Tcheby. Furthermore, in order to
verify the generality of the modified Tchebycheff
aggregation function, MOEA/D-DE is combined with the proposed aggregation which is called MOEA/D-DE-tmd. The original MOEA/D-DE is compared with MOEA/D-DE-tmd in the comparative experiments.

\begin{figure*}
\centering
  \includegraphics[width=6.5in]{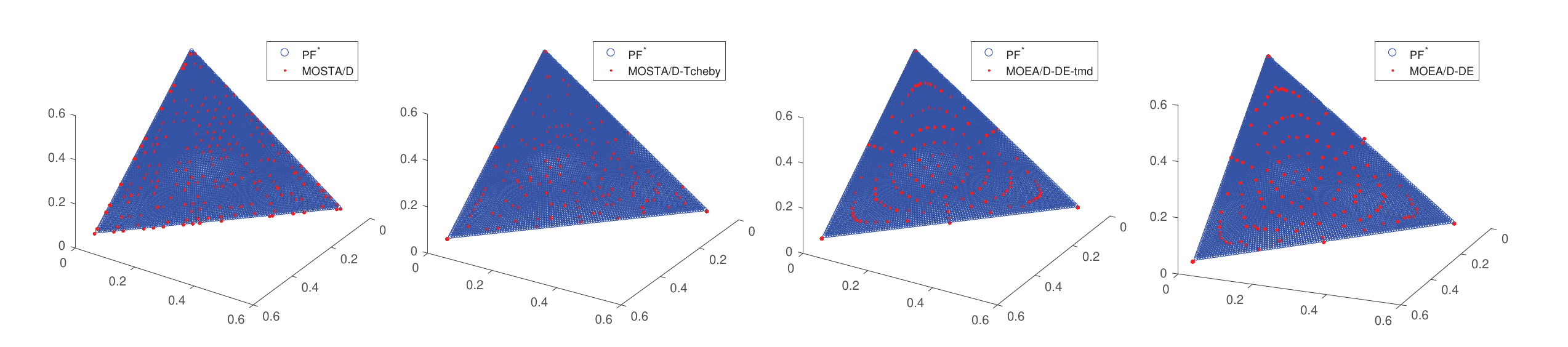}\\
  \caption{The Pareto fronts of P1 test function obtained by different algorithms}
   \label{DTLZ1}
\end{figure*}
\begin{figure*}
\centering
  \includegraphics[width=6.5in]{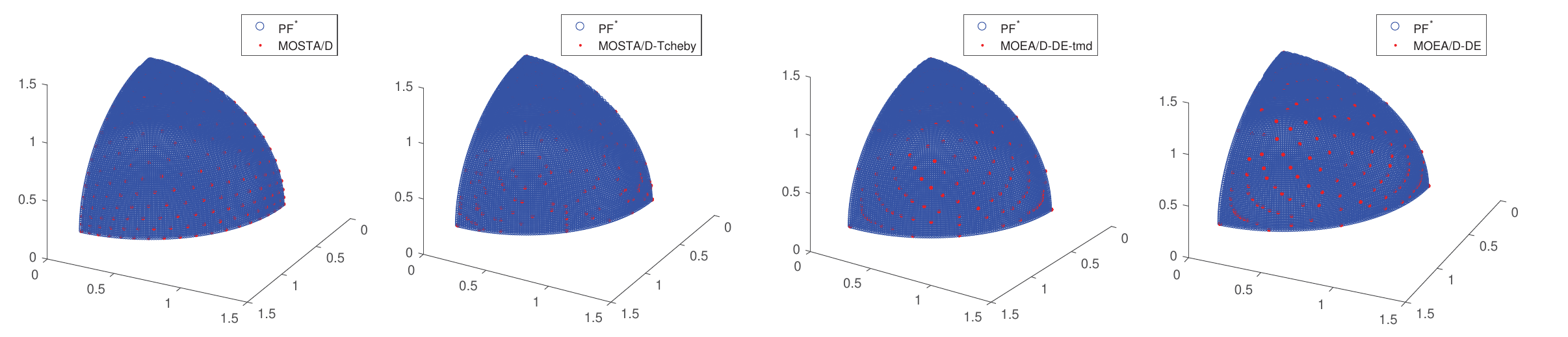}\\
  \caption{The Pareto fronts of P2 test function obtained by different algorithms}
   \label{DTLZ2}
\end{figure*}
\begin{figure*}
\centering
  \includegraphics[width=6.5in]{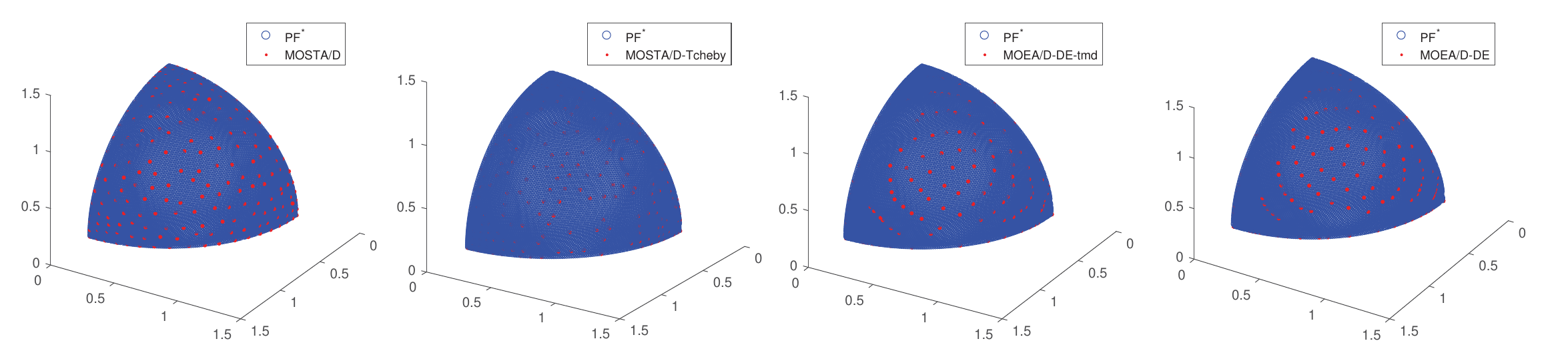}\\
  \caption{The Pareto fronts of P3 test function obtained by different algorithms}
  \label{DTLZ4}
\end{figure*}
\begin{figure*}
\centering
  \includegraphics[width=6.5in]{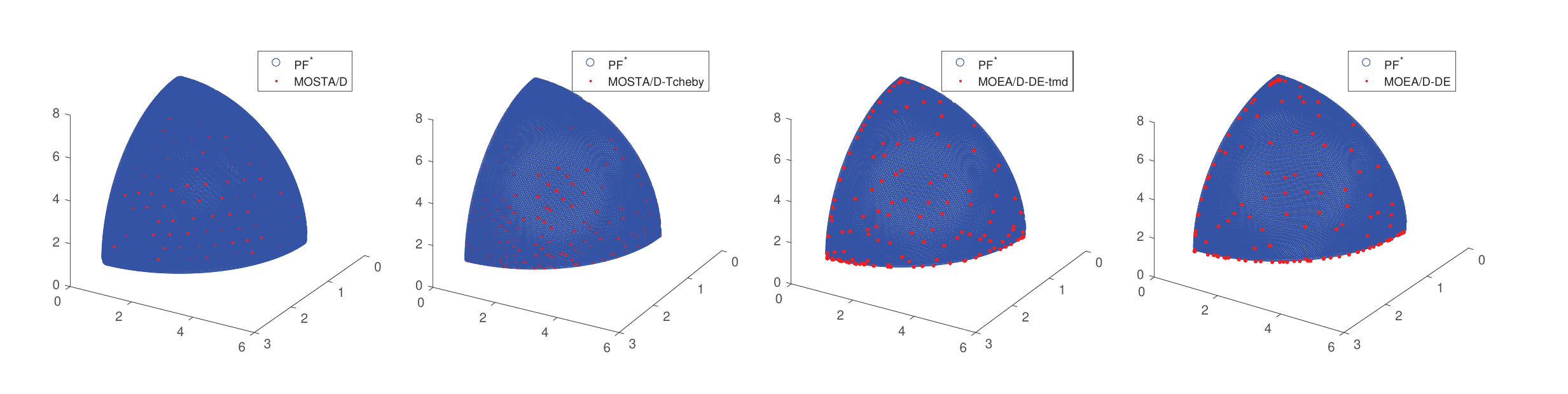}\\
  \caption{The Pareto fronts of P4 test function obtained by different algorithms}
\label{WFG4}
\end{figure*}

\begin{table*}
  \centering
  \caption{Statistical results of the IGD$^+$ values obtained by four algorithms}
  \label{IGD-verify}
\scalebox{0.8}{\begin{tabular}{ccccccccccccccc}
\toprule
problem & MOSTA/D & MOSTA/D-Tcheby  & MOEA/D-DE-tmd    & MOEA/D-DE   \\
\midrule
\multirow{1}*{P1} & \textbf{1.2499e-2 (1.4161e-1)} &  1.3566e-2 (1.5612e-1)$-$ & 1.4622e-2 (4.6521e-1)$-$ & 1.5062e-2 (3.6250e-4)$-$\\
\multirow{1}*{P2} & \textbf{1.7191e-2 (5.9085e-5)} &  1.7651e-2 (6.5443e-5)$-$ & 2.1156e-2 (6.7319e-4)$-$& 2.4279e-2 (2.4264e-4)$-$\\
\multirow{1}*{P3} &\textbf{ 1.0067e-2 (2.0978e-4)} &  1.0099e-2 (5.6215e-4)$\approx$ & 1.1936e-2 (1.5414e-3)$-$ & 1.2048e-2 (1.5201e-3)$-$\\
\multirow{1}*{P4} & \textbf{ 1.4201e-1 (1.3655e-4) }&  1.4295e-1 (6.3129e-4)$-$ & 1.4416e-1 (2.3641e-3)$-$ & 1.4647e-1 (3.9685e-3)$-$\\
\bottomrule
\end{tabular}}
\end{table*}

\begin{table*}
  \centering
  \caption{Statistical results of the HV values obtained by four algorithms}
  \label{HV-verify}
\scalebox{0.8}{\begin{tabular}{ccccccccccccccc}
\toprule
problem& MOSTA/D & MOSTA/D-Tcheby  & MOEA/D-DE-tmd    & MOEA/D-DE   \\
\midrule
\multirow{1}*{P1} &\textbf{ 8.4224e-1 (8.9115e-4)} &  8.3511e-1 (7.6611e-4)$-$ & 8.3356e-1 (2.3654e-4)$-$ & 8.2671e-1 (6.4312e-4)$-$\\
\multirow{1}*{P2} & \textbf{5.7603e-1 (8.0312e-5) }& 5.6225e-1 (6.4652e-5)$-$ & 5.6125e-1 (9.3525e-4)$-$ & 5.5197e-1 (8.3625e-4)$-$\\
\multirow{1}*{P3} & \textbf{5.5568e-1 (7.9064e-2) }& 5.5535e-1 (5.3245e-2)$\approx$  & 5.5477e-1 (6.5565e-2)$-$ & 5.5436e-1 (2.4214e-2)$-$\\
\multirow{1}*{P4} & \textbf{5.3329e-1 (1.1525e-3) }&  5.3036e-1 (4.5512e-3)$-$ & 5.2595e-1 (6.9462e-3)$-$ & 5.0595e-1 (2.5334e-3)$-$\\
\bottomrule
\end{tabular}}
\end{table*}

Figs.\ref{DTLZ1}-\ref{WFG4} plots the distribution of the final solutions with the median IGD$^+$-metric value of all algorithms for each test instance. Obviously, four algorithms can converge to the true PFs of P1-P4.
In addition, MOSTA/D can approximate the PFs of these problems quite well.

The statistical results  are summarized in Table \ref{IGD-verify} and Table \ref{HV-verify}, respectively.  MOSTA/D is much better than comparison algorithms on all instances. For P1, MOSTA/D performs better than its variants in terms of the mean of IGD$^+$ and HV values. MOEA/D-DE is outperformed by MOEA/D-DE-tmd. For P2, the performance of MOSTA/D and MOEA/D-DE-tmd on this instance is much better than MOSTA/D-Tcheby and MOEA/D-DE, respectively.
By contrast, MOSTA/D-Tcheby is slightly outperformed by MOSTA/D on P3.
P4 is a test problem with multimodality. MOSTA/D obtains a higher HV value and a smaller IGD$^+$ value than other three algorithms on P4. Furthermore, the performance of MOEA/D-DE-tmd is improved compared to MOEA/D-DE.

From the above discussion, MOSTA/D achieves the best results and shows the most competitive overall performance on all instances. MOEA/D-DE combined with the modified Tchebycheff aggregation function is still significantly better than the original MOEA/D-DE. The improved algorithms both achieve better performance because the proposed aggregation function comprehensively takes into account the Tchebycheff decomposition approach and the relationship.  The better candidate solutions are retained by considering the matching degree during the selection process. It is evident that the proposed aggregation function shows adequate generality and outperforms the traditional Tchebycheff aggregation significantly. MOSTA/D can generate evenly distributed solutions.

\subsection{Validation of the proposed algorithm}
In this part, several benchmark test functions and a typical
engineering optimization problem are adopted to verify the
performance of the proposed algorithm.

\subsubsection{Benchmark test function verification}
10 benchmark test functions (P5-P14) are adopted to verify the proposed algorithm compared with MOEA/DD  \cite{li2014evolutionary}, MOEA/D-DE \cite{li2009multiobjective}, MOEA/D-M2M
\cite{liu2013decomposition} and NSLS \cite{chen2014new}.

\begin{figure*}

\centering
  \includegraphics[width=6.5 in]{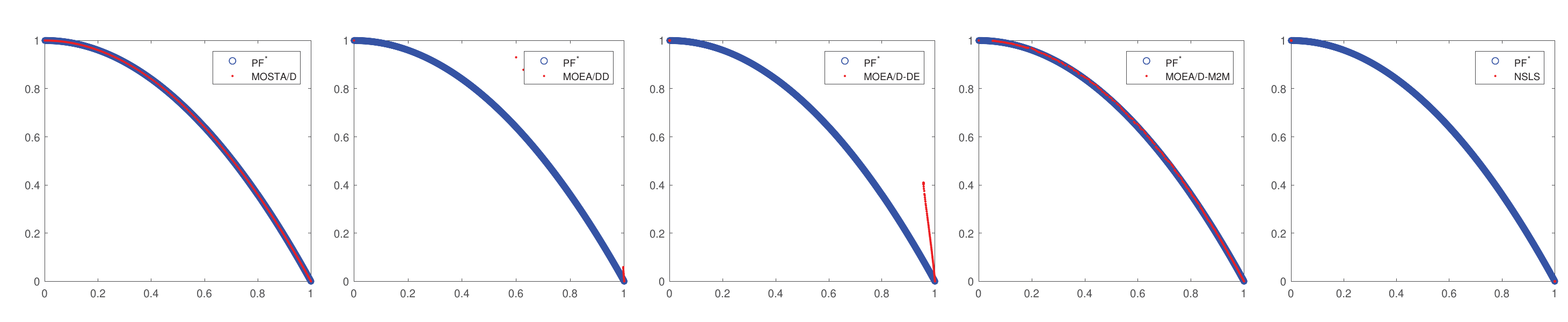}\\
  \caption{The Pareto fronts of P5 test function obtained by different algorithms}
  \label{MOP2}
\end{figure*}

\begin{figure*}
\centering
  \includegraphics[width=6.5in]{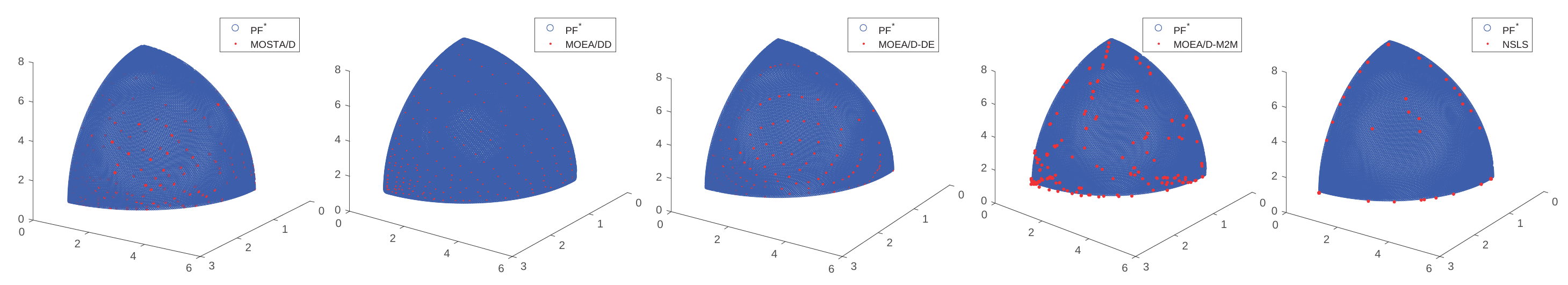}\\
  \caption{The Pareto fronts of P14 test function obtained by different algorithms}
  \label{WFG6}
\end{figure*}

The PFs are  plotted to visualize the Pareto optimal solutions obtained by algorithms in Figs.\ref{MOP2}-\ref{WFG6}.  The statistical results obtained by five algorithms are summarized in Table \ref{IGD-test} and Table \ref{HV-test}, respectively. The PFs of both P5 and P6 are non-convex, and MOSTA/D has
achieved the best performance  on these two problems. Obviously, MOSTA/D can find a set of evenly distributed solutions. However, MOEA/DD, MOEA/D-DE, and NSLS can only obtain boundary points.  P7 is a relatively complex test problem, only MOSTA/D and MOEA/D-M2M can generate evenly distributed solutions on P7. These two algorithms can obtain all the branches of the PF, while other algorithms completely fail to obtain the true PF.

MOSTA/D is tried on P8-P10 with more decision variables. Although the PF of P8 is disconnected, MOSTA/D can solve this problem with the minimal IGD$^+$-metric and the maximal HV-metric, which means  MOSTA/D  has good distribution and convergence. MOEA/DD, MOEA/D-DE and MOEA/D-M2M can also solve this problem. However, it appears that NSLS can only find boundary points of the true PF.
It can be observed that the overall performance of  each algorithm  is generally good. Although MOEA/D-M2M is slightly outperformed by MOSTA/D, its performance is still significantly better than MOEA/D-DE and NSLS.
MOEA/D-DE and NSLS can not guarantee that all solutions converge to the optimal solution on P10. It can be observed that MOSTA/D succeeds in experiments for more decision variables wtih good scalability.

The PF approximated by MOSTA/D is of high quality on a highly multimodal problem P11, although the performance of MOEA/DE, MOEA/D-M2M, and NSLS are not very stable on this problem,
as evidenced by the results, since all of them completely can not reach the true PF. MOSTA/D is better than other comparison algorithms on P12 with a disconnected PF. However, the other three algorithms can not generate evenly distributed solutions on this problem. MOSTA/D is slightly outperformed by MOEA/DD on P13. Compared to MOEA/D-DE, MOEA/D-M2M, and NSLS, the performance of MOSTA/D on this instance is much better. MOSTA/D has achieved a comparable smaller IGD$^+$ value and higher HV value than comparison algorithms on P14. MOEA/DD  also show generally competitive performance.

Since different state transformation operators in MOSTA/D can be used for global search, local search and heuristic search. Alternative use of local search and global search can help MOSTA/D quickly converge to the PFs. By considering the matching relationship between weight vectors and candidate solutions, the proposed algorithm can obtain better convergence and diverse solution set than comparsion algorithms.

\begin{table*}[!htb]
  \centering
  \caption{Statistical results  of the IGD$^+$  values obtained by five algorithms}
  \label{IGD-test}
\scalebox{0.8}{\begin{tabular}{cccccccccccc}
\toprule
                    problem  & MOSTA/D  & MOEA/DD  & MOEA/D-DE    & MOEA/D-M2M  & NSLS   \\
\midrule
\multirow{1}*{P5} &\textbf{ 1.1481e-3 (2.2140e-7)} & 1.4564e-1 (5.3836e-3)$-$ & 1.3025e-1 (1.5715e-2)$-$  & 7.5914e-2 (3.1244e-2)$-$ & 1.4909e-1 (3.7785e-3)$-$ \\
\midrule
\multirow{1}*{P6} & \textbf{1.0619e-3 (6.0176e-7)} & 1.8225e-1 (2.6033e-2)$-$ & 2.1911e-1 (3.6865e-2)$-$  & 7.9932e-2 (3.7211e-2)$-$  & 9.9584e-2 (2.8230e-17)$-$ \\
\midrule
\multirow{1}*{P7} & \textbf{7.7437e-4  (2.6381e-6)} & 1.5455e-1 (3.2050e-2)$-$ & 2.0563e-1 (9.2561e-3)$-$ & 1.7721e-2 (1.9433e-2)$-$ & 2.3283e-1 (9.1057e-3)$-$ \\
\midrule
\multirow{1}*{P8} & \textbf{9.9740e-4 (2.8634e-5)} & 1.6625e-2 (3.1993e-3)$-$ & 2.0534e-2 (1.0725e-1)$-$ & 1.1956e-3 (6.6263e-5)$-$ & 5.1385e-2 (2.4295e-2)$-$\\
\midrule
\multirow{1}*{P9} & \textbf{1.0604e-3 ( 8.5456e-6)} & 4.9521e-3 (8.5396e-4)$-$ & 1.3045e-3 (5.3767e-5)$-$ & 1.0969e-3 (1.2366e-4)$\approx$  & 3.1531e-3 (1.9472e-3)$-$\\
\midrule
\multirow{1}*{P10} & \textbf{1.0715e-2 (4.5265e-3)} & 6.7252e-2 (2.0754e-3)$-$ & 2.1526e-1 (2.1124e-2)$-$ & 2.5237e-1 (2.8539e-2)$-$ &  6.6245e-1 (5.1078e-2)$-$\\
\midrule
\multirow{1}*{P11} & \textbf{4.6781e-1 (1.1551)} & 9.0413e-1 (1.0963)$-$ & 1.1752 (3.6271)$-$ &  23.1653 (14.3258)$-$ & 161.5125 (10.4172)$-$\\
\midrule
\multirow{1}*{P12} & \textbf{2.9161e-2 (1.0541e-3)} & 2.9419e-2 (1.1764e-3) $-$ & 4.9734e-2 (1.4172e-3)$-$  & 9.1199e-2 (4.7485e-3)$-$ & 5.6519e-1 (3.2921e-2)$-$ \\
\midrule
\multirow{1}*{P13} & 1.2301e-1 (2.6418e-4) & \textbf{1.2280e-1 (9.9305e-4)}$+$ & 1.6187e-1 (1.3274e-3)$-$ & 1.6416e-1 (3.7584e-3)$-$ & 5.3582e-1 (5.9314e-2)$-$\\
\midrule
\multirow{1}*{P14} & \textbf{1.3500e-1 (1.6452e-2)} & 1.3582e-1 (2.1244e-2)$\approx$ & 1.4646e-1 (7.2599e-2)$-$ & 3.1364e-1 (8.6691e-2)$-$ & 3.4252e-1 (1.0409e-2)$-$ \\
\bottomrule
\end{tabular}}
\end{table*}

\begin{table*}[!htb]
  \centering
  \caption{Statistical results of the HV values obtained by five algorithms}
  \label{HV-test}
\scalebox{0.8}{\begin{tabular}{cccccccccccc}
\toprule
                    problem  & MOSTA/D  & MOEA/DD  & MOEA/D-DE    & MOEA/D-M2M  & NSLS   \\
\midrule
\multirow{1}*{P5} & \textbf{7.7096e-1 (1.2901e-5)} & 3.2329e-1 (1.3325e-2)$-$ & 2.2690e-1 (3.6212e-2)$-$ & 3.2206e-1 (5.5134e-2)$-$  & 1.7870e-1 (8.0923e-3)$-$ \\
\midrule
\multirow{1}*{P6} &\textbf{ 6.5248e-1 (1.7623e-5)} & 9.5423e-2 (6.5641e-3)$-$ & 9.0909e-2 (7.0622e-17)$-$ & 2.3000e-1 (5.5312e-2)$-$ & 1.7355e-1 (2.8234e-17)$-$ \\
\midrule
\multirow{1}*{P7} & \textbf{ 9.5648e-1 (1.6514e-5)}& 3.7145e-1 (3.4156e-2)$-$ & 3.0932e-1 (1.6943e-2)$-$ & 5.7433e-1 (2.7222e-2)$-$  & 2.5948e-1 (1.4821e-2)$-$ \\
\midrule
\multirow{1}*{P8} & \textbf{1.1198 (9.5667e-5)} & 1.0639 (2.7493e-1)$-$ & 1.0906 (1.6071e-1)$-$ & 1.1085 (2.0602e-2)$\approx$ & 1.0226 (4.7033e-2)$-$ \\
\midrule
\multirow{1}*{P9} & \textbf{1.3579 (2.1326e-4)}& 9.7559e-1 (4.5909e-2)$-$  & 1.3575 (1.3452e-4)$\approx$ & 1.3571 (2.6147e-2)$-$  & 1.3560 (1.7251e-3)$-$\\
\midrule
\multirow{1}*{P10} &\textbf{ 85.9517 (2.2931e-1)} & 85.2654 (7.1665e-1)$-$  & 84.3821 (1.9492e-1)$-$  & 85.4620 (5.3654e-1)$-$ & 82.4864 (7.6543e-1)$-$ \\
\midrule
\multirow{1}*{P11} & \textbf{5.7326e-1 (1.7415e-3) }& 5.7066e-1 (1.9384e-1)$-$  & 5.1106e-1 (1.3943e-1)$-$  & 0.0000e+0 (0.0000e+0)$-$  & 0.0000e+0 (0.0000e+0)$-$ \\
\midrule
\multirow{1}*{P12} & \textbf{9.3688e-1 (1.2416e-3)} & 9.3671e-1 (1.0532e-3)$\approx$  & 9.0723e-1 (4.1623e-3)$-$  & 9.0658e-1 (4.2724e-3)$-$   & 6.6931e-1 (1.8134e-2)$-$ \\
\midrule
\multirow{1}*{P13} & 5.2220e-1 (1.9412e-3) & \textbf{5.2622e-1 (6.4660e-4)}$+$ & 4.8639e-1 (1.3632e-3)$-$ &4.9167e-1 (3.1134e-3)$-$ & 3.6388e-1 (9.3423e-3)$-$\\
\midrule
\multirow{1}*{P14} & \textbf{5.1965e-1 (2.0525e-2)} & 5.1892e-1 (1.4363e-2)$-$ & 4.9008e-1 (6.0913e-2)$-$ & 3.9609e-1 (4.8223e-2)$-$ & 3.8638e-1 (7.1432e-3)$-$\\
\bottomrule
\end{tabular}}
\end{table*}

\begin{figure*}
\centering
  \includegraphics[width=6.5in]{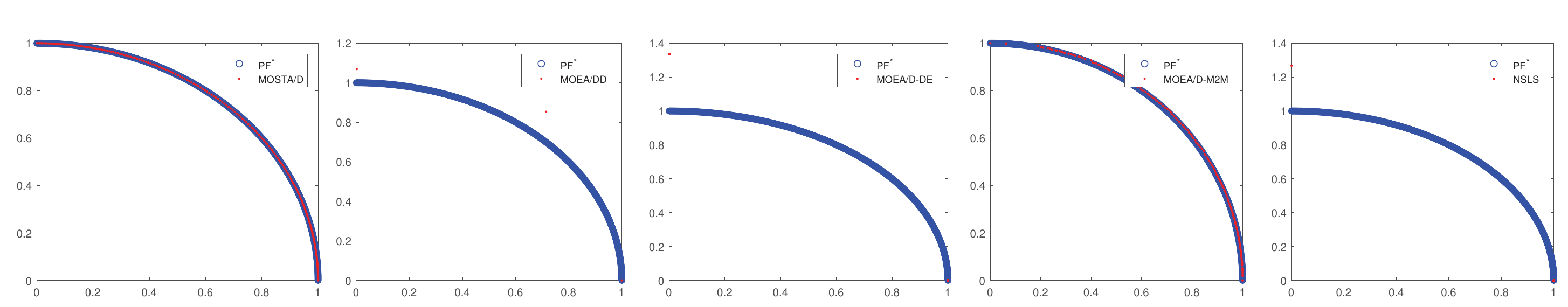}\\
  \caption{The Pareto fronts of P6 test function obtained by different algorithms}
\end{figure*}

\begin{figure*}
\centering
  \includegraphics[width=6.5in]{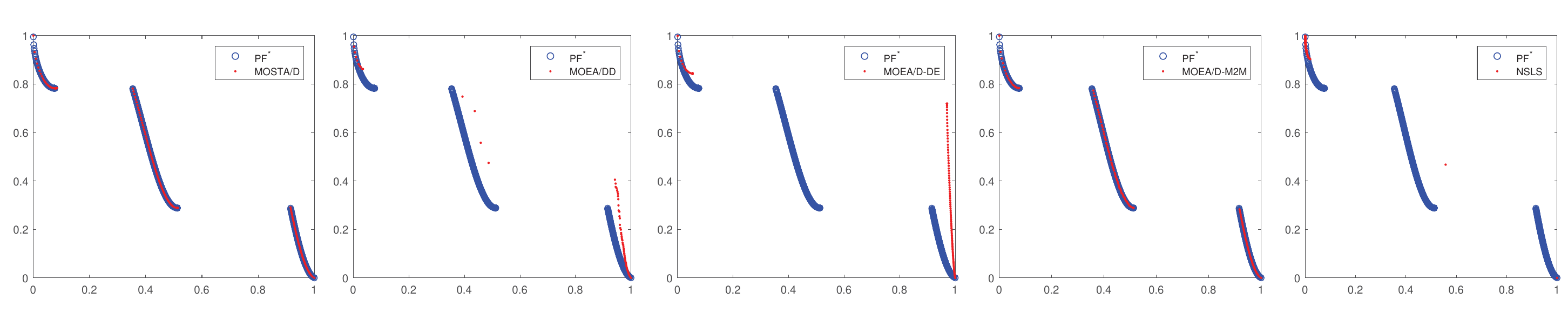}\\
  \caption{The Pareto fronts of P7 test function obtained by different algorithms}
\end{figure*}

\begin{figure*}
\centering
  \includegraphics[width=6.5in]{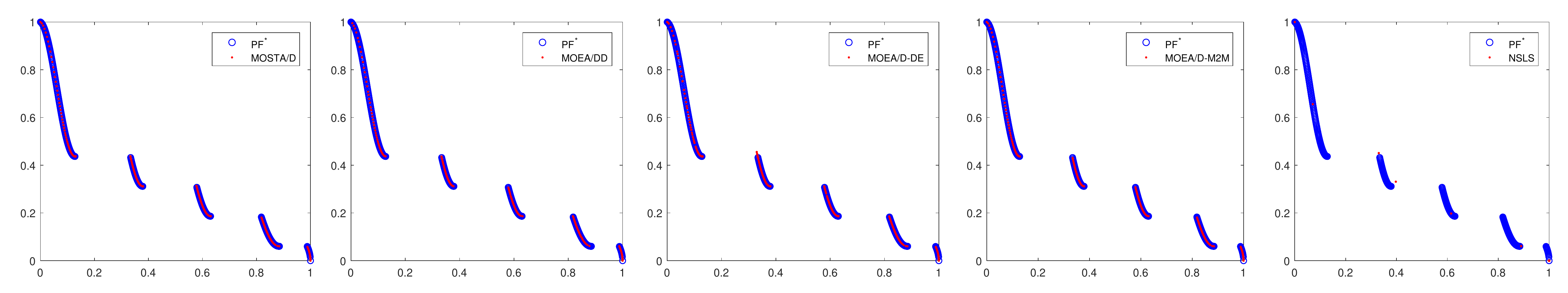}\\
  \caption{The Pareto fronts of P8 test function obtained by different algorithms}
\end{figure*}

\begin{figure*}
\centering
  \includegraphics[width=6.5in]{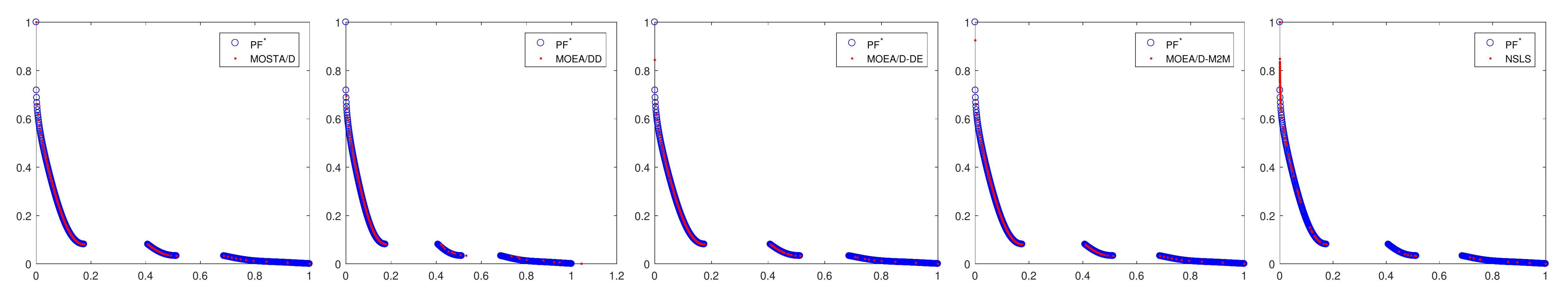}\\
  \caption{The Pareto fronts of P9 test function obtained by different algorithms}
\end{figure*}

\begin{figure*}
\centering
  \includegraphics[width=6.5in]{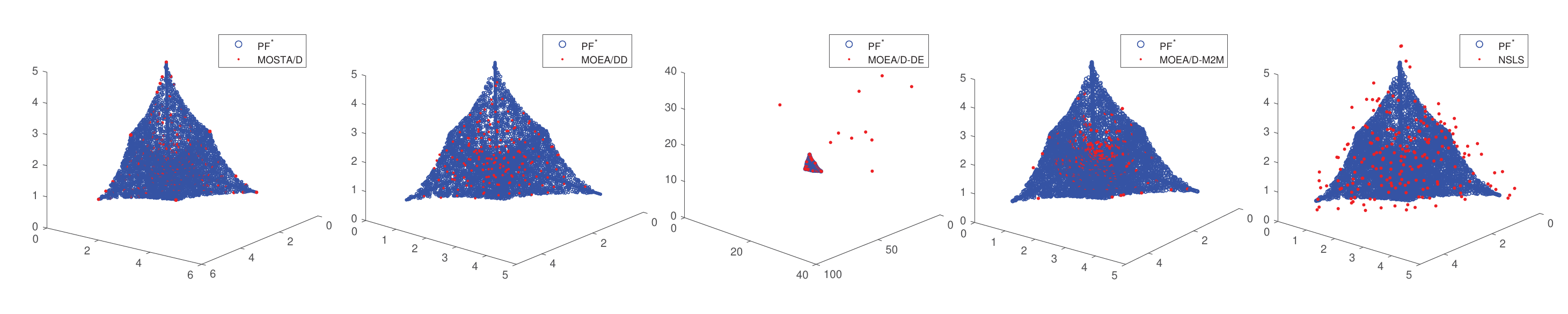}\\
  \caption{The Pareto fronts of P10 test function obtained by different algorithms}
\end{figure*}
\begin{figure*}
\centering
  \includegraphics[width=6.5in]{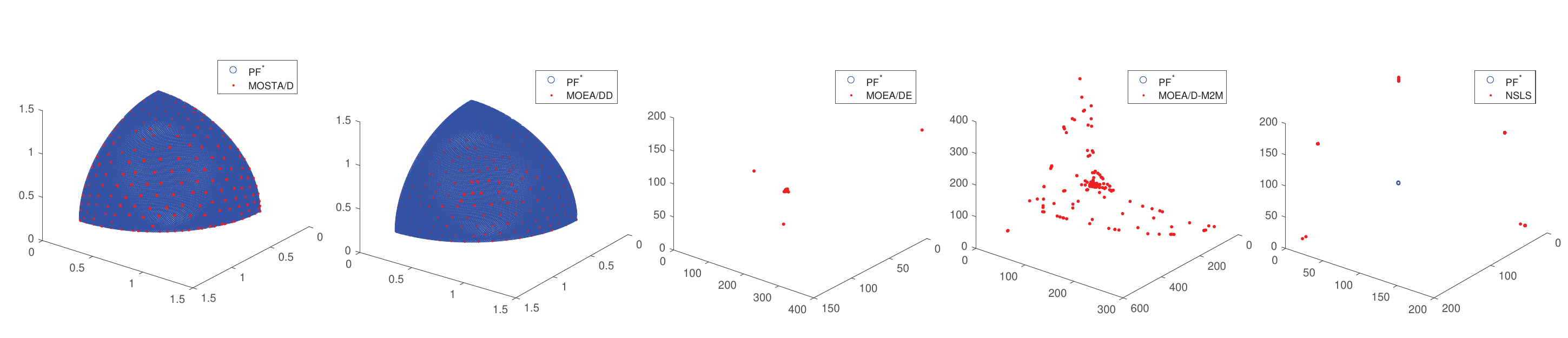}\\
  \caption{The Pareto fronts of P11 test function obtained by different algorithms}
\end{figure*}

\begin{figure*}
\centering
  \includegraphics[width=6.5in]{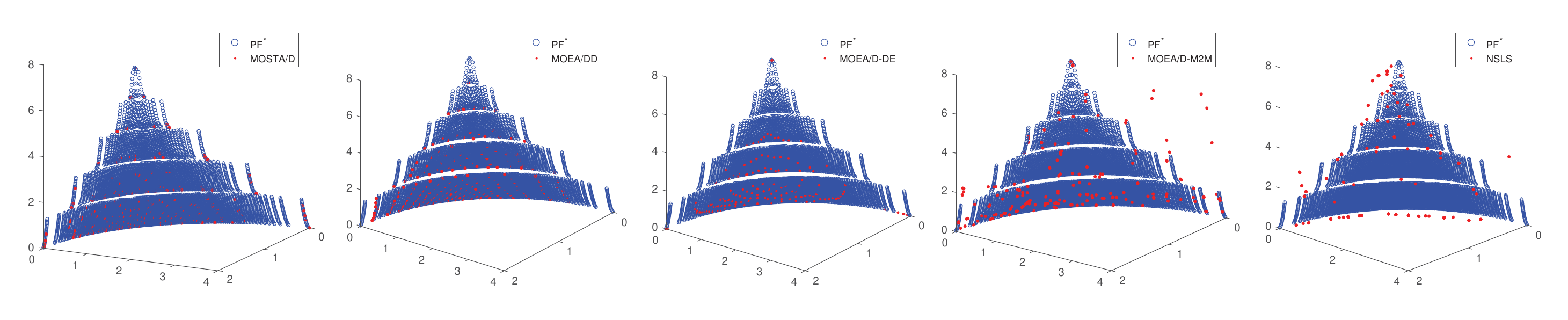}\\
  \caption{The Pareto fronts of P12 test function obtained by different algorithms}
\end{figure*}

\begin{figure*}
\centering
  \includegraphics[width=6.5in]{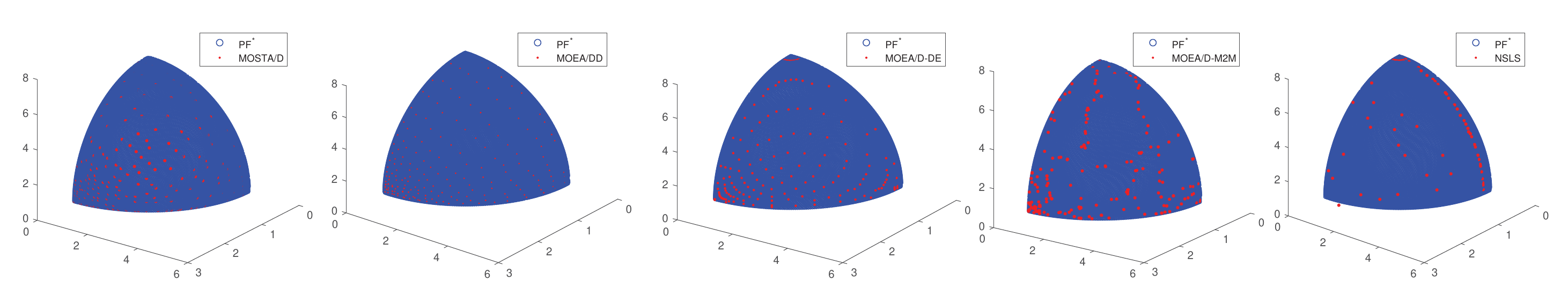}\\
  \caption{The Pareto fronts of P13 test function obtained by different algorithms}
\end{figure*}

\subsubsection{Engineering optimization problem verification}
\begin{figure}
\centering
  \includegraphics[width=5.5cm]{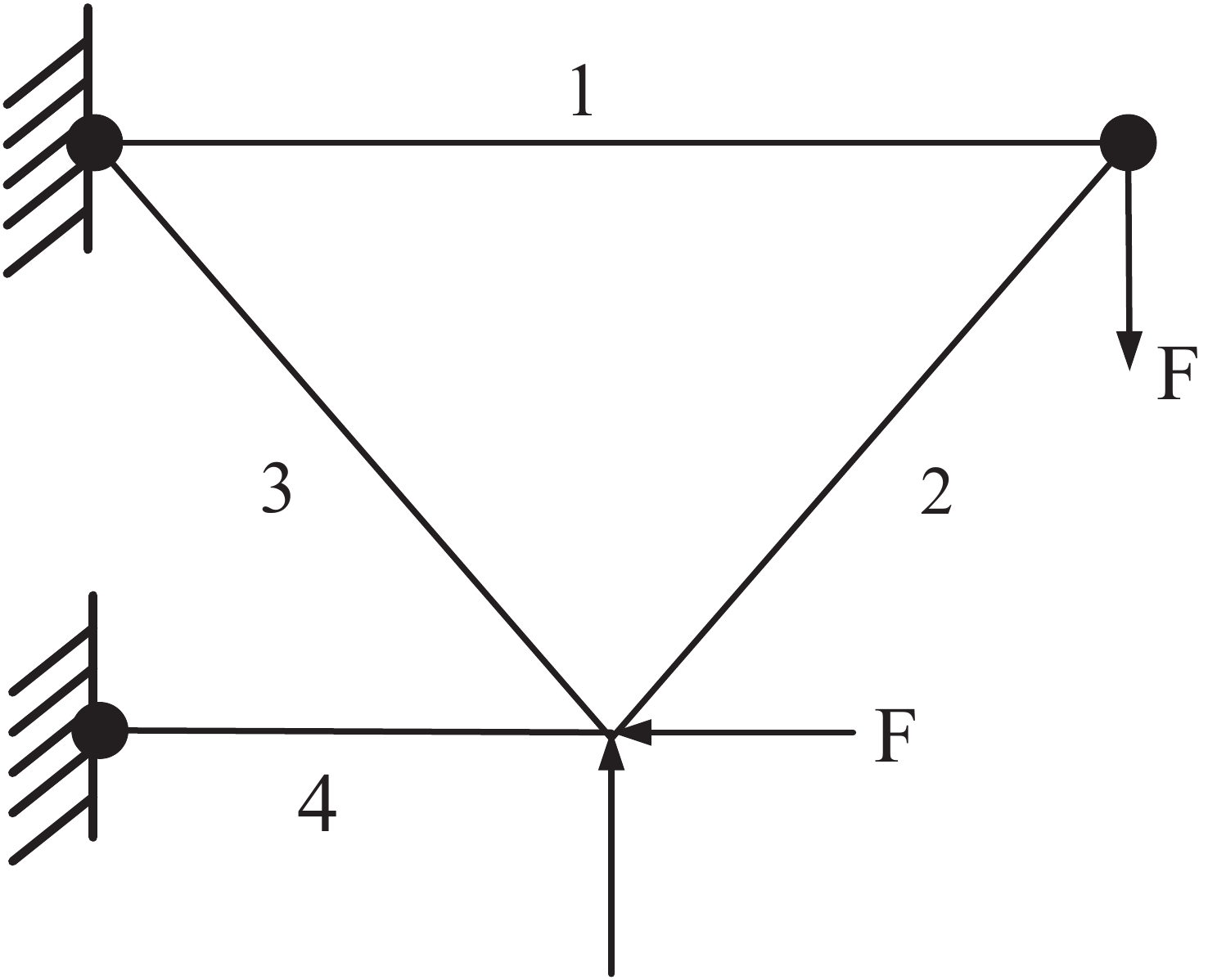}\\
  \caption{The diagrammatic sketch of four bar plane truss}
  \label{fourbar1}
\end{figure}

The four bar plane truss design is a typical optimization problem in the structural optimization field,
in which structural mass ($f_1$) and compliance ($f_2$) of a 4-bar plane truss should be minimized.
The design drawing are shown in Fig. \ref{fourbar1}.
The four bar plane truss design optimization problem can be formulated
as follows:

\begin{eqnarray}
\label{fourbar}
&&f_1(x)=L(2x_1+\sqrt{2}x_2+\sqrt{x_3}+x_4)\\ \nonumber
&&f_2(x)=\frac{FL}{E}(\frac{2}{x_1}+\frac{2\sqrt{2}}{x_2}- \frac{2\sqrt{2}}{x_3}+\frac{2}{x_4}) \\ \nonumber
&&\mathrm{s.t.}  \quad (F / \epsilon  )\leq x_1 \leq 3(F / \epsilon  )\\ \nonumber
&& \quad \quad \sqrt{2}(F / \epsilon  ) \leq x_2 \leq 3(F / \epsilon  )\\ \nonumber
&& \quad \quad\sqrt{2}(F / \epsilon ) \leq x_3 \leq 3(F / \epsilon  )\\ \nonumber
&& \quad \quad(F / \epsilon  ) \leq x_4 \leq 3(F / \epsilon  )\\ \nonumber
\end{eqnarray}
where $F=10kN$, $E = 2*10^5kN/cm^2$, $\epsilon  = 10kN/cm^2$, $L = 0.2m$.

\begin{figure*}
\centering
  \includegraphics[width=6in]{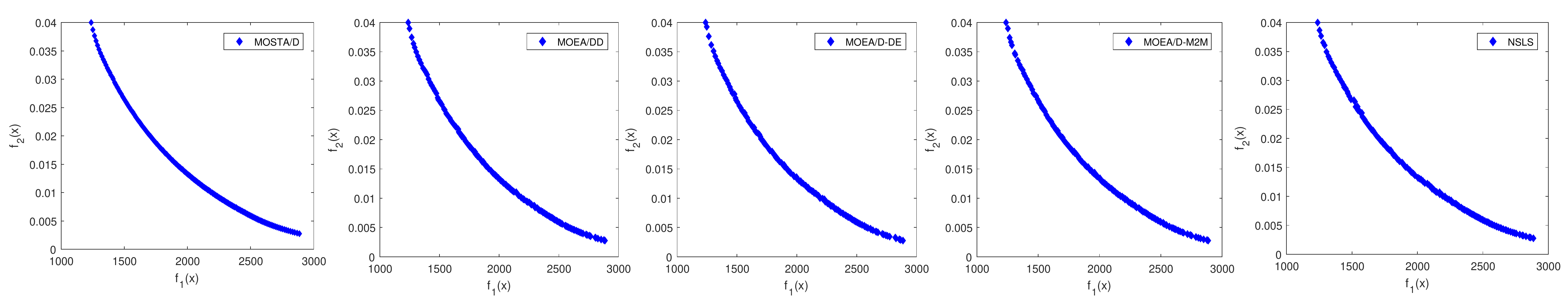}\\
  \caption{The Pareto fronts of four bar plane truss design problem obtained by different algorithms}
  \label{fourbarParetosnew}
\end{figure*}

It can be observed that the objective functions differ greatly in order of magnitude. Hence, a parameter is added to the second objective function and the optimization objectives can be shown as Eq. (\ref{fourbar2}). After optimization, the objective function value is transformed into the original optimization objective function value.

\begin{eqnarray}
\label{fourbar2}
&&f_1(x)=L(2x_1+\sqrt{2}x_2+\sqrt{x_3}+x_4)\\ \nonumber
&&f_2(x)=\frac{\tau FL}{E}(\frac{2}{x_1}+\frac{2\sqrt{2}}{x_2}- \frac{2\sqrt{2}}{x_3}+\frac{2}{x_4}) \\ \nonumber
\end{eqnarray}
where $\tau=75000$.

The PFs formed by the set of target values corresponding to Pareto optimal solution sets obtained by five algorithms are shown in Fig. \ref{fourbarParetosnew}.
Compared with the PF in \cite{chiandussi2012comparison} and four other MOEAs, the PF obtained by MOSTA/D  has relatively better diversity and convergence. Therefore, it can be said that the proposed MOSTA/D is able to solve practical engineering optimization problems. In the future, the proposed algorithm will be applied to solve more complex practical problems in industry, including multi-objective optimization-based PID control, $H_\infty$ and robust control of process control. How to design a reasonable and effective weight generation strategy will be further studied in the future for solving higher dimensional MOPs.

\section{Conclusion}
In this paper, the influence of matching relationship between
weight vectors and  candidate solutions on selection and updating procedure in decomposition approaches for solving MOPs
are studied and analyzed. Considering that the matching
relationship has great influence and it is difficult to quantify,
the concept of matching degree is proposed based on vectorial
angle among the reference point vector, candidate solutions
vectors and weight vectors. Based on the matching degree, a new
decomposition approach is proposed which comprehensively
takes into account the Tchebycheff approach and matching
degree. Moreover, it can be proved that the proposed decomposition approach has the same functionality with the Tchebycheff approach.
Furthermore, a decomposition based multi-objective state
transition algorithm is proposed. By testing several benchmark test functions and an engineering optimization problem, MOSTA/D can approximate the true PF
with high convergence precision and good diversity.


%

%

\ifCLASSOPTIONcompsoc
  \section*{Acknowledgments}
\else
  \section*{Acknowledgment}
\fi

This study was funded by the National Natural Science Foundation of China (Grant No. 61873285),
International Cooperation and Exchange of the National Natural Science Foundation of China (Grant No. 61860206014) and the National Key Research and Development Program of China (Grant No. 2018AAA0101603).

\ifCLASSOPTIONcaptionsoff
  \newpage
\fi

\normalem

\bibliographystyle{IEEEtran}
\bibliography{mybibfile}

\end{document}